\documentclass[12pt]{amsart}

\usepackage{amsmath, amssymb, amscd}
\usepackage[mathscr]{eucal}
\usepackage{mathrsfs}
\usepackage{verbatim}
\usepackage{version}
\usepackage{nicefrac}
\usepackage{pstricks}
\usepackage{pst-all}
\usepackage{color}
\usepackage[all]{xy}
\usepackage{mathdots}
\usepackage{longtable}

\newenvironment{NB}{
\color{blue}{\bf NB}. \footnotesize }{}
\excludeversion{NB}

\newcommand{\bM}{\mathbf M}
\newcommand{\bN}{\mathbf N}

\newcommand{\bP}{\mathbf P}

\newcommand{\bS}{\mathbf S}
\newcommand{\bU}{\mathbf U}

\newcommand{\cG}{\mathcal G}

\newcommand{\cK}{\mathcal K}

\newcommand{\cM}{\mathcal M}
\newcommand{\cN}{\mathcal N}
\newcommand{\cC}{\mathcal C}
\newcommand{\cO}{\mathcal O}

\newcommand{\cX}{\mathcal X}

\newcommand{\fg}{\mathfrak g}

\newcommand{\fo}{\mathfrak{o}}

\newcommand{\ft}{\mathfrak{t} }

\newcommand{\fsp}{\mathfrak{sp} }

\newcommand{\heta}{\hat {\eta}}

\newcommand{\hR}{\hat R}

\newcommand{\oast}{{\overline{\ast}}}

\newcommand{\og}{\overline{g}}

\newcommand{\oPhi}{{\overline{\Phi}}}

\newcommand{\obS}{\overline{\bS}}

\newcommand{\obU}{\overline{\bU}}

\newcommand{\rO}{\mathrm{O}}

\newcommand{\te}{{\widetilde e}}

\newcommand{\tf}{{\widetilde f}}

\newcommand{\tM}{\tilde M}
\newcommand{\tm}{\tilde m}

\newcommand{\tp}{{\widetilde p}}

\newcommand{\tS}{\tilde S}

\newcommand{\tX}{\tilde X}

\newcommand{\aaa}{\mathbb A}
\newcommand{\pp}{\mathbb P}

\newcommand{\cc}{\mathbb C}
\newcommand{\rr}{\mathbb R}
\newcommand{\zz}{\mathbb Z}

\newcommand{\ep}{\epsilon}
\newcommand{\vep}{\varepsilon}

\newcommand{\ad}{\mathrm{ad} }

\newcommand{\diag}{\mathrm{diag} }

\newcommand{\dual}{^\vee }

\newcommand{\End}{\mathrm{End} }

\newcommand{\gl}{\mathfrak{gl} }
\newcommand{\gV}{\mathfrak{g}(V) }

\newcommand{\GL}{\mathrm{GL} }

\newcommand{\git}{/\!\!/ }
\newcommand{\Hom}{\mathrm{Hom} }

\newcommand{\Ker}{\mathrm{Ker} }
\newcommand{\Image}{\mathrm{Im} }

\newcommand{\Lie}{\mathrm{Lie} }
\newcommand{\Mat}{\mathrm{Mat} }
\newcommand{\inv}{^{-1}}
\newcommand{\Id}{\mathrm{Id} }
\newcommand{\OY}{\mathrm{O}.Y}
\newcommand{\fp}{\mathfrak{p} }

\newcommand{\rank}{\mathrm{rank} }
\newcommand{\reg}{\mathrm{reg} }

\newcommand{\sm}{\mathrm{sm} }
\newcommand{\pt}{\mathrm{pt}}

\newcommand{\SL}{\mathrm{SL}}

\newcommand{\Sp}{\mathrm{Sp}}
\newcommand{\Spin}{\mathrm{Spin}}
\newcommand{\SO}{\mathrm{SO}}
\newcommand{\SU}{\mathrm{SU}}
\newcommand{\tr}{\mathrm{tr} }
\newcommand{\USp}{\mathrm{USp}}

\newcommand{\I}{{\mathrm{I}}}
\newcommand{\II}{{\mathrm{II}}}
\newcommand{\III}{{\mathrm{III}}}
\newcommand{\IV}{{\mathrm{IV}}}

\newtheorem{prop}{Proposition}[section]
\newtheorem{thm}[prop]{Theorem}
\newtheorem{lem}[prop]{Lemma}

\newtheorem{cor}[prop]{Corollary}

\theoremstyle{remark}
  \newtheorem{rk}[prop]{Remark}
  
\theoremstyle{definition}
 
 \newtheorem{defn}[prop]{Definition}

\numberwithin{equation}{section}

\begin{document}


\title[Moduli spaces of framed symplectic and orthogonal bundles]{Moduli spaces of framed symplectic and orthogonal bundles on $\pp^2$ and the K-theoretic Nekrasov partition functions}
\author{Jaeyoo Choy}
\address{Dept. Math., Kyungpook Nat'l Univ., Sangyuk-dong, Buk-gu, Daegu 702-701, Korea}
\email{choy@knu.ac.kr}

\dedicatory{Dedicated to the late Professor Kentaro Nagao}


\begin{abstract}
Let $K$ be the compact Lie group $\USp(N/2)$ or $\SO(N,\rr)$. Let $\cM_n^K$ be the moduli space of framed $K$-instantons over $S^4$ with the instanton number $n$. By \cite{Do}, $\cM_n^K$ is endowed with a natural scheme structure. It is a Zariski open subset of a GIT quotient of $\mu\inv(0)$, where $\mu$ is a holomorphic moment map such that $\mu\inv(0)$ consists of the ADHM data.

The purpose of the paper is to study the geometric properties of $\mu\inv(0)$ and its GIT quotient, such as complete intersection, irreducibility, reducedness and normality. If $K=\USp(N/2)$ then $\mu$ is flat and $\mu\inv(0)$ is an irreducible normal variety for any $n$ and even $N$. If $K=\SO(N,\rr)$ the similar results are proven for low $n$ and $N$.

As an application one can obtain a mathematical interpretation of the K-theoretic Nekrasov partition function of \cite{NS}.
\end{abstract}

\maketitle
\setcounter{tocdepth}{1}
\tableofcontents

\thispagestyle{empty} \markboth{Jaeyoo Choy}{Moduli spaces of symplectic and orthogonal vector bundles and the K-theoretic Nekrasov partition function}


\section{Introduction}
\label{sec: intro}


\subsection{}
\label{subsec: intro: different motivation}

The Nekrasov partition function was formulated by Nekrasov \cite{Nek} in the 4-dimensional $\cN=2$ supersymmetric gauge theory in physics, especially relevant to the Seiberg-Witten prepotential \cite{SW}. It is defined as a generating function of the equivariant integration of the trivial cohomology class 1 over Uhlenbeck partial compactification of the (framed) instanton moduli space on $\rr^4$ for all the instanton numbers.  Its logarithm turned out to contain the Seiberg-Witten prepotential as the coefficient of the lowest degree (with respect to two variables from the torus action on $\rr^4$), which is the remarkable result proven by Nakajima-Yoshioka \cite{NY}, Nekrasov-Okounkov \cite{NO} (both for $\SU(N)$) and Braverman-Etingof \cite{BE} (for any gauge group) in completely independent methods.
%
%
The Nekrasov partition function is an equivariant version of Donaldson-type invariants for $\rr^4$, where the ordinary Donaldson invariants are integrals over Uhlenbeck compactifications of instanton moduli spaces on compact 4-manifolds. More generally, a Nekrasov partition function for the theory with matters is defined via the equivariant integration of a cohomology class other than the trivial class. For instance, in the works of G\"ottsche-Nakajima-Yoshioka \cite{GNY1}\cite{GNY3}, such partition functions are used to express wallcrossing terms of Donaldson invariants and a relation between them and Seiberg-Witten invariants.

Our interest in this paper is the K-theoretic Nekrasov partition function. It arose from the 5-dimensional $\cN=2$ supersymmetric gauge theory proposed by Nekrasov \cite{Nek96}. Mathematically, one change in its definition is made \cite{Nek}: the equivariant integration of elements in the equivariant K-theory of coherent sheaves instead of cohomology classes. There is a technical, but a subtle problem here. We need a scheme structure on the Uhlenbeck space for the definition, but there are several choices (see \cite{BE}). Therefore it is not clear what is the correct definition. For type A, one can use framed moduli spaces of coherent sheaves instead, which are smooth. The K-theoretic Nekrasov partition functions also appeared in mathematics literatures as K-theoretic Donaldson invariants studied in \cite{NYII}\cite{GNY2} in which the K-theory version of Nekrasov conjecture and blowup equations are proven for type A.

Nekrasov-Shadchin \cite{NS} took a different approach when the gauge group is of classical type. They defined the partition function using the K-theory class of the Koszul complex which defines the ADHM data of the instanton moduli spaces given in \cite{ADHM}.

The main purpose of this paper is to show that Nekrasov-Shadchin's definition coincides with a generating function of the coordinate rings of the instanton moduli spaces for the classical gauge groups. The answer for the gauge group $\SU(N)$ has been already known by a general result of Crawley-Boevey on quiver varieties \cite{CB2}. (See \S\ref{subsec: motivating questions} for the precise motivating questions which we pursue). Our study for the gauge groups $\USp(N/2)$ and $\SO(N,\rr)$ informs of some geometry of the moduli spaces in algebraic geometry, as the instanton moduli space for $\USp(N/2)$ and $\SO(N,\rr)$ is isomorphic to the moduli space of (framed) vector bundles with symplectic and orthogonal structures on $\pp^2$ respectively. (See \S\ref{intro: subsec: ADHM2} for explanation on the scheme structures of moduli spaces.)


\subsection{}
\label{subsec: intro AHS}

Let us fix the notation and explain earlier results in order to state our result precisely. Let $K$ be a compact connected simple Lie group. Let $P_K$ be a principal $K$-bundle over the 4-sphere $S^4$. Since $\pi_3(K)\cong\zz$, an integer $n$ uniquely determines the topological type of $P_K$. Let $\cM_n^K$ be the quotient of the space of ASD connections (instantons) by the group of the gauge transformations trivial at infinity $\infty\in S^4$ \cite[\S5.1.1]{DK}. By \cite[Table 8.1]{AHS} $\cM^K_n$ is a $\cC^\infty$-manifold with $\dim \cM^K_n=4nh^\vee$, where $h^\vee$ is is the dual Coxeter number of $K$. In this paper, we will also consider the case $K=\SO(4,\rr)$. Its universal cover $\Spin(4)$ is the product $\SU(2)\times \SU(2)$, and hence $\SO(4,\rr)$-bundles over $S^4$ are classified by pairs of integers $(n_L,n_R)$. We have $\cM^K_{(n_L,n_R)}\cong \cM^{\SU(2)}_{n_L}\times\cM^{\SU(2)}_{n_R}$, and hence $\dim \cM^K_{(n_L,n_R)}=8(n_L+n_R)$.


\subsection{}\label{intro: subsec: ADHM}

Let $K=\SU(N)$. Donaldson \cite{Do} showed that $\cM^K_n$ is naturally isomorphic to the moduli space of framed holomorphic vector bundles $E$ with $c_2(E)=n$ over $\pp^2$, where a framing is a trivialization of $E$ over the line $l_\infty$ at infinity. This was proved by using the ADHM description \cite{ADHM} and the relation between moment maps and geometric invariant theory (GIT).
Let $K=\USp(N/2)$ or $\SO(N,\rr)$ and $G$ be its complexification. Since $G$-bundles can be identified with rank $N$ vector bundles with symplectic or orthogonal structures, we have the following description of $\cM^K_n$.

Let $G_k':=\rO(k)$ (resp.\ $\Sp(k/2)$) if $G=\Sp(N/2)$ (resp.\ $\SO(N)$). Let $\fp(\cc^k)$ be the space of symmetric endomorphisms of $\cc^k$ (see \eqref{eq: t p} for the precise definition).
The ADHM data are elements $x\in \bN$ satisfying $\mu(x)=0$, where   $\bN:=\fp(\cc^k)^{\oplus2}\oplus \Hom(\cc^N,\cc^k)$ and $\mu$ is the holomorphic moment map defined on $\bN$. The $G_k'$-action on $\bN$ induces a $G_k'$-action on $\mu\inv(0)$. Now $\cM^K_n$ is the image of the regular locus $\mu\inv(0)^\reg$ of the GIT quotient $\mu\inv(0)\to \mu\inv(0)\git G_k'$.
Here the second Chern number $k$ of associated vector bundles is determined from $n$ by the argument in \cite[\S10]{AHS} as
    \begin{equation}\label{eq: rel of c2 and instanton number}
    k=\left\{ \begin{array}{lll} n & & \mbox{if $G=\Sp(N/2)$,}
        \\
        2n & & \mbox{if $G=\SO(N)$, where $N\ge5$,}
        \\
        4n & & \mbox{if $G=\SO(3)$.}
        \end{array}\right.
    \end{equation}

The GIT quotient $\mu\inv(0)\git G_{k}'$ endows the Donaldson-Uhlenbeck partial compactification of $\cM^{K}_{n}$ with a scheme structure except the case $G=\SO(3)$ (cf.\ \S\ref{subsec: stratification}).
We call $\mu\inv(0)\git G_{k}'$ the scheme-theoretic Donaldson-Uhlenbeck partial compactification of $\cM^{K}_{n}$ unless $G=\SO(3)$.
As we will see in \S\ref{subsec: motivating questions}, the K-theoretic Nekrasov partition function coincides with the generating function of the coordinate rings of algebraic functions on these Donaldson-Uhlenbeck spaces if $G=\Sp(N/2)$.


\subsection{}\label{intro: subsec: ADHM2}

The scheme structure of $\cM^K_n$ is given by $\mu\inv(0)^\reg/G_{k}'$. Since $\mu\inv(0)^\reg/G_{k}'$ is a Zariski-open subset of $\mu\inv(0)\git G_{k}'$, $\cM^K_n$ is a quasi-affine scheme.
We notice that $\cM^K_n$ may have two scheme structures because it is possible that $\Lie(K)=\Lie(K')$ for a different classical group $K'$ (hence $\cM^K_n=\cM^{K'}_n$) but $\cM^K_n$ and $\cM^{K'}_n$ have the different ADHM descriptions. Such pairs of $(K,K')$ are $(\SU(2),\USp(1))$, $(\SU(2),\SO(3,\rr))$,  $(\USp(2),\SO(5,\rr))$ and $(\SU(4),\SO(6,\rr))$. A scheme-theoretic isomorphism $\cM^K_n\stackrel\cong\to \cM^{K'}_n$ is induced by sending the associated vector bundles $E$ to itself, $\ad E$, $(\Lambda^2 E)_0$ and $\Lambda^2 E$ respectively. The notations in the above are as follows: $\ad E$ is the trace-free part of $\End(E)$, and
$(\Lambda^2 E)_0$ is the kernel of the symplectic form $\Lambda^2E\to \cO$.

In the above isomorphism we used the following two assertions: First by Donaldson's theorem \cite{Do}, $\cM^K_n$ is canonically isomorphic to the moduli space of framed rank $N$ vector bundles $E$ with $c_2(E)=k$ and $G$-structure. 
Then the latter space is endowed with a (smooth) scheme structure and the canonical isomorphism is scheme-theoretic. See \S\ref{sec: moduli space of G bundles}.


\subsection{}\label{intro: subsec: group action}

There is a canonical $G\times G_k'$-action on $\bN$.
The canonical $G$-action on $\bN$ induces a $G$-action on $\cM^K_n$. We define a $(\cc^*)^2$-action on $\bN$ by $(q_1,q_2).(B_1,B_2,i)=(q_1B_1,q_2B_2,\sqrt{q_1q_2}i)$. More precisely this action is well-defined only on a double cover of $\cc^*\times \cc^*$, but we follow the convention in physics. It commutes with the $G\times G_k'$-action and induces a $(\cc^*)^2$-action on $\cM^K_n$. Since $i$ always appears together with $i^*$ for generators in $\cc[\bN]^{G_k'}$ (see \cite[Thereoms 2.9.A and 6.1.A]{Weyl}), the $(\cc^*)^2$-action is well-defined on $\cM^K_n$.

Let $T_G$ be a maximal torus of $G$. We identify $T_G=(\cc^*)^l$, where $l=\rank(G)$. Let $T:=T_G\times (\cc^*)^2=(\cc^*)^{l+2}$. Then $T$ acts on $\mu\inv(0)\git G_k'$ and $\cM^K_n$.
Nekrasov-Shadchin \cite{NS} defined the $K$-theoretic instanton partition
function $Z^K$ for the gauge group $K$ by an explicit integral
formula. Their formula can be interpreted as follows.

Fix the instanton number $n$.
Let us consider $\mu$ as a section of the trivial vector bundle $\mathbf
E = \operatorname{Lie}G'_k\times\mathbf N$ over $\mathbf N$. We endow
$\mathbf E$ with a $G'_k\times T$-equivariant structure so that $\mu$
is a $T$-equivariant section. Then $\mu$ defines a Koszul complex
\begin{equation}\label{eq:Koszul}
    \Lambda^{\operatorname{rank} \mathbf E} \mathbf E^\vee
    \to \cdots \to
    \Lambda^2 \mathbf E^\vee \to
    \mathbf E^\vee \to
    \Lambda^0\mathbf E^\vee = \mathcal O_{\mathbf N}.
\end{equation}
The alternating sum
\(
   \sum_i (-1)^i \Lambda^i\mathbf E^\vee
\)
of terms gives an element in $K^{T\times G'_k}(\mathbf N)$, the Grothendieck group of $T\times G'_k$-equivariant vector bundles over $\mathbf N$.
We then take a pushforward of the class with respect to the obvious
map $p\colon \mathbf N\to \mathrm{pt}$. This is not a class in the
representation ring $R({T\times G'_k})$ because $\mathbf N$ is not
proper.
However, it is a well-defined class in $\operatorname{Frac}(R(T\times G'_k))$, the fractional field of $R(T\times G'_k)$ by {\it equivariant
  integration}: Check first that the origin $0$ is the unique fixed point
in $\mathbf N$ with respect to $T\times T'_k$, where $T'_k$ is a maximal
torus of $G'_k$. (See \S\ref{sec: finite dim}\begin{NB}
 Replace by the apropriate label
\end{NB}.)
Therefore the pushforward homomorphism $p_*$ can be defined as the
inverse of $i_*$ thanks to the fixed point theorem of the equivariant
K-theory, where $i\colon \{0\}\to\mathbf N$ is the inclusion. In
practice, we take the Koszul resolution of the skyscraper sheaf at the
origin as above, replacing $\mathbf E$ by $\mathbf N$, and then divide
the class by $\sum_i (-1)^i\Lambda^i \mathbf N^\vee$.

In our situation, $\mathbf N$ has only positive weights with respect
to the factor $(\mathbb C^*)^2$ of $T$. (See \S\ref{sec: finite dim}\begin{NB}
 Replace by the apropriate label
\end{NB}.)
Therefore the rational function can be expanded as an element in a
{\it completed\/} ring of $T\times G'_k$-characters
\begin{equation*}
    \hat R(T\times G'_k) := R(T_G\times G'_k)[[q_1^{-1},q_2^{-1}]],
\end{equation*}
where $q_m$ denote the $T$-characters of the $1$-dimensional
representations for the second factor $(\mathbb C^*)^2$ of $T = T_G\times (\mathbb C^*)^2$.
The expansion is considered as a formal character, as the sum of
dimensions of weight spaces: Each weight space is finite dimensional
and its dimension is given by the coefficient of $f\in\hat R(T\times
G'_k)$ of the monomial corresponding to the weight.

Next we take the $G'_k$-invariant part. This is given by an integral
over the maximal compact subgroup of $T'_k$ by Weyl's integration
formula.
Finally taking the generating function for all $n\in\mathbb Z_{\ge 0}$
with the formal variable $\mathfrak q$, we define the instanton
partition function
\begin{equation*}
    Z^K := \sum_{n\ge 0} \mathfrak q^n
    \sum_i (-1)^i p_*(\Lambda^i \mathbf E^\vee)^{G'_k}
    \in \hat R({T})[[\mathfrak q]],
\end{equation*}
where $\hat R(T) := R(T_G)[[q_1^{-1},q_2^{-1}]]$.

Explicit computation of $Z^K$ has been performed for small instanton
numbers in physics literature. See \cite{BHM}\cite{HKS} for instance. \begin{NB} Replace by the
    apropriate label. Also I delete 13, 19 as they do not contain
    relevant computation. Add 1012.4468 instead.
\end{NB}


\subsection{}\label{subsec: motivating questions}

Nekrasov-Shadchin's definition of $Z^K$ is mathematically rigorous,
but depends on the ADHM description, hence is not intrinsic in the
instanton moduli space $\mathcal M^K_n$. We can remedy this problem
under the following geometric assumptions:
\begin{enumerate}
      \item $\mu$ is a flat morphism;
      \item $\mu^{-1}(0)\git G'_k$ is a normal variety and $\mathcal
    M^K_n$ has complement of codimension $\ge 2$ in
    $\mu^{-1}(0)\git G'_k$.
\end{enumerate}
Under these assumptions, $\sum_i (-1)^i p_*(\Lambda^i \mathbf
E^\vee)^{G'_k}$ is the formal character of the ring $\mathbb
C[\mathcal M^K_n]$ of regular functions on $\mathcal M^K_n$:
First, (1) asserts $\mu$ gives a regular sequence so that the Koszul
complex \eqref{eq:Koszul} is a resolution of the coordinate ring
$\mathbb C[\mu^{-1}(0)]$ of $\mu^{-1}(0)$. Secondly (2) asserts
$\mathbb C[\mathcal M_n^K] = \mathbb C[\mu^{-1}(0)\git G'_k] =
\mathbb C[\mu^{-1}(0)]^{G'_k}$ by extension of regular functions on
$\mathcal M^K_n = \mu^{-1}(0)^{\operatorname{reg}}/G'_k$ to
$\mu^{-1}(0)\git G'_k$.
As explained in \S\ref{intro: subsec: ADHM2}, the scheme structure of $\cM^K_n$ is independent of the choice of the ADHM description.

When $K=\mathrm{SU}(N)$, (1),(2) are proved by Crawley-Boevey \cite{CB} in a general context of quiver varieties.
Moreover the Gieseker space $\widetilde{\mathcal M}^{N}_n$ of framed rank $N$ torsion free sheaves $\mathcal E$ on $\mathbb P^2$ with $c_2(\mathcal E)=n$ gives a resolution of singularities of $\mu^{-1}(0)\git G'_k$. 
The higher direct image sheaves $R^{l}\tp_{*}\cO_{\widetilde{\mathcal M}^{N}_n}=0$ for $l\ge1$, where $\tp\colon \
\widetilde{\mathcal M}^{N}_n\to \pt$ is the obvious map.
For, $\widetilde{\mathcal M}^{N}_n$ is symplectic and thus the canonical sheaf $\cK_{\widetilde{\mathcal M}^{N}_n}$ is isomorphic to $\cO_{\widetilde{\mathcal M}^{N}_n}$.
By the Grauert-Riemenschneider vanishing theorem, $R^{l}\pi_{*}\cK_{\widetilde{\mathcal M}^{N}_n}=0$ for $l\ge1$, where $\pi$ denotes the resolution of singularities.
Now the vanishing of higher direct image sheaves comes from the $E_{2}$-degeneration of the spectral sequence $R^{l}p_{*}R^{m}\pi_{*} \Rightarrow R^{l+m}\tp_{*}$, where $p\colon \mu^{-1}(0)\git G'_k\to \pt$ is the obvious affine map.
Therefore $\mathbb C[\widetilde{\mathcal M}^{N}_n]
\cong \mathbb C[\mu^{-1}(0)\git G'_k]$.
The $T$-action lifts to $\widetilde{\mathcal M}^{N}_n$, where fixed
points are para\-metrized by $N$-tuples of Young diagrams corresponding
to direct sums of monomial ideal sheaves.
Therefore the character of $\mathbb C[\widetilde{\mathcal M}^{N}_n]$
is given by a sum over $N$-tuples of Young diagrams, which is the original definition of the instanton partition function in \cite{Nek}. \begin{NB}
    In the current version, it is 34.
\end{NB}
See \cite{NYII} for more detail on the latter half of this argument. \begin{NB}
    In the current version, it is 32.
\end{NB}


\subsection{}

A goal of this paper is to prove (1),(2) for
$K=\mathrm{USp}(N/2)$. (See Theorem \ref{th: main 5}.\begin{NB}
 Replace by the apropriate label
\end{NB}) A key of the proof is a result of Panyushev \cite{Pa}, which gives the flatness of $\mu$ for $N=0$.

We also study $\mu^{-1}(0)$ for $K = \mathrm{SO}(N,\mathbb R)$ for
$(N,k) = (2,k)$, $(N,2)$ or $(3,4)$. (See Theorems \ref{th: N2}, \ref{th: k2}, \ref{th: N3k4}
respectively.\begin{NB}
 Replace by the apropriate label
\end{NB})
The case $N=2$ is less interesting since there are no
$\mathrm{SO}(2,\mathbb R)$-instantons except for $k=0$. But
$\mu^{-1}(0)$ {\it does\/} make sense, hence the study of its
properties is a mathematically meaningful question. We show that $\mu$
is not flat, i.e., (1) is not true in this case.
Similarly there is no instanton for $(N,k) = (3,2)$. We will see that
$\mu$ is flat, but $\mu^{-1}(0)\git G'_k$ with reduced scheme structure is isomorphic to $\mathbb C^2$, so (2) is false in this case.
In the case $(N,2)$ with $N\ge 4$, we prove that $\mu^{-1}(0)\git
G'_k$ is isomorphic to the product of $\mathbb C^2$ and the closure
$\mathbf P$ of the minimal nilpotent $\mathrm{O}(N)$-orbit.
When $(N,k) = (3,4)$, we show that $\mu$ is flat and
$\mu^{-1}(0)\git G'_k$ is a union of two copies of $\mathbb
C^2\times\mathbf P$ meeting along $\mathbb C^2\times \{0\}$. Here
$\mathbf P$ is the minimal (= regular) nilpotent $\mathrm{O}(3)$-orbit
closure, which is isomorphic to $\mathbb C^2/\mathbb Z_2$. In
particular, (2) is false.
The isomorphism $\mathcal M^{\mathrm{SU}(2)}_n\cong\mathcal
M^{\mathrm{SO}(3)}_n$ implies that both $\mathcal M^{\mathrm{SU}(2)}_n$ and $\mathcal
M^{\mathrm{SO}(3)}_n$ are $\mathbb
C^2\times(\mathbf P\setminus\{0\})$ for $n=1$.
\begin{NB}
    replace an appropriate statement, if you can say about the scheme
    structure.
\end{NB}%
However their ADHM
descriptions are different, as it is $\mathbb C^2\times\mathbf P$ for
$\mathrm{SU}(2)$. This phenomenon happens only for
$\mathrm{SO}(3)$. See Theorem \ref{th: main 4}. \begin{NB}
 Replace by the apropriate label
\end{NB}
These examples show that the definition of $Z^K$
depends on the ADHM description, and hence must be studied with care.

The author plans to study the case $(N\ge 4, k\ge 4)$ in the near future.


\subsection{}
\label{subsec: intro contents}

The paper is organized as follows.

In \S\ref{set up} we give the set-up for the entire part of the paper. Specifically we set up linear algebra of vector spaces with nondegenerate forms in \S\ref{subsec: right adj}--\S\ref{subsec: **}.
In \S\ref{subsec: the Orbit-closedness and the semisimplicity} we identify the closed $G_k'$-orbits in $\bN$. By this identification we see that $\cM^K_n$ is the $G_k'$-orbit space of stable-costable representations in $\mu\inv(0)$.
In \S\ref{subsec: stratification}
we stratify $\mu\inv(0)\git G_k'$ with
$\cM^K_n$ as a stratum. Each stratum is isomorphic to a product of $\cM^K_{n'}$ and a symmetric product of $\cc^2$ for some $n'\le n$.

In \S\ref{sec: geometry of the moduli spaces of framed symplectic bundles} and \S\ref{sec: geometry of the moduli spaces of framed orthogonal bundles} we state Theorems \ref{th: main 5}, \ref{th: k2}--\ref{th: N3k4} which describes the geometry of $\mu\inv(0)$ when $K=\USp(N/2)$ and $\SO(N,\rr)$.

In \S\ref{subsec: Kraft-Procesi's classification theory of nilpotent pairs} we explain Kraft-Procesi's classification theory of nilpotent pairs. In our study of $\mu\inv(0)$, their theory is quite useful to see contribution from the factor $\Hom(\cc^N,\cc^k)$ of $\bN$. In the case $k=2$, the geometry of $\mu\inv(0)$ is immediately deduced from Kraft-Procesi's theory on $\Hom(\cc^N,\cc^k)$ since $\fp(\cc^2)$ consists of the scalars.

In \S\ref{subsec: proof Theorem k2} and \S\ref{subsec: proof Theorem N2} we will see how Kraft-Procesi's theory is applied to the cases $k=2$ and $N=2$. This will prove Theorem \ref{th: k2} and a part of Theorem \ref{th: N2}.

In \S\ref{sec: the moduli spaces of framed orthogonal bundles of rank 3 and k4} we prove Theorem \ref{th: N3k4} (the case $(N,k)=(3,4)$). The proof involves more than Kraft-Procesi's theory since the pairs $(B_1,B_2)\in \fp(\cc^4)^{\oplus2}$ are no more commuting pairs. Since $[B_1,B_2]\neq0$ in general, the study of $\mu\inv(0)$ does not solely come from the factor $\Hom(\cc^3,\cc^4)$. So we study the commutator map $\fp(\cc^4)^{\oplus2}\to \Lie(\Sp(2))$, $(B_1,B_2)\mapsto [B_1,B_2]$.

In \S\ref{sec: proof of Theorem N2 (2)} we finish the proof of Theorem \ref{th: N2} (the case $(N,k)=(2,4)$) based on the study of the commutator map in the above.

\vskip.3cm \noindent\textit{Acknowledgements.}
This work was supported by JSPS RONPAKU (Dissertation Ph.D.) Program. I would like to express my deepest gratitude to my supervisor Professor Hiraku Nakajima. His guidance, care, patience and enthusiasm have been leading me to pursue this topic and to learn many neighboring subjects. I thank Research Institute for Mathematical Sciences, Kyoto University, for nice environment for research during my visits over years. I also thank Professor Bumsig Kim for hospitality at Korea Institute for Advanced Study where a part of the paper is written. I thank Professors Yuji Tachikawa, Ugo Bruzzo and Matt Young for their interest in this work.
I thank the anonymous referee for careful reading and suggestions which led to this revision of the paper. 


\section{Preliminary}
\label{set up}

We set up conventions and notation and review basic material for an entire part of the paper.

We are working over $\cc$.
Vector spaces are all finite dimensional and schemes are of finite type. We say that a morphism between schemes is irreducible (resp.\ normal and Cohen-Macaulay) if all the nonempty fibres are irreducible (resp.\ normal and Cohen-Macaulay). If $\cM$ is a scheme then $\cM^{\mathrm{sm}}$ and $\cM^{\mathrm{sing}}$ are the smooth locus and the singular locus of $\cM$ respectively.

Let $G$ be an algebraic group and $\fg:=\Lie(G)$ the Lie algebra of $G$. Let $\cM$ be a $G$-scheme. Let $G^x:=\{g\in G|\, g.x=x\}$ the stabilizer subgroup of $x\in \cM$. Let $\fg^x:=\Lie(G^x).$


\subsection{The right adjoint}
\label{subsec: right adj}
Let $V_1$ and $V_2$ be vector spaces with nondegenerate bilinear forms $(\,,\,)_{V_1}$ and $(\,,\,)_{V_2}$ respectively. Then for any $i\in \Hom(V_1,V_2)$, we have the right adjoint $i^*\in \Hom(V_2,V_1)$, i.e. $(v,i^*w)_{V_1}=(iv,w)_{V_2}$, where $v\in V_1$ and $w\in V_2$. The map
    $$
    *:\Hom(V_1,V_2)\to \Hom(V_2,V_1),\   i\mapsto i^*
    $$
is a $\cc$-linear isomorphism. Further if $V_3$ is a vector space with a nondegenerate bilinear form then for $i\in \Hom(V_1,V_2),\ j\in \Hom(V_2,V_3)$, we have $(ji)^*=i^*j^*$.


\subsection{Anti-symmetric and symmetric forms}
\label{subsec: as s}

Let $V$ be a vector space of dimension $k$ with a nondegenerate form $(\,,\,)_V$. Let $\vep\in \{-1,+1\}$. Let $(\ ,\ )_\vep$ be a nondegenerate bilinear form $(u,v)_\vep=\vep(v,u)_\vep$ on $V$. If $\vep=+1$ (resp.\ $-1$) then $(\ ,\ )_\vep$ is an orthogonal form (resp.\ symplectic form). We say $V$ is orthogonal (resp.\ symplectic) if $\vep=+1$ (resp.\ $\vep=-1$).

We decompose $\gl(V)=\ft(V)\oplus \fp(V)$ as a vector space, where
    \begin{equation}
    \label{eq: t p}
    \begin{aligned}
    \ft(V) &:= \{X\in \gl(V)|\, (Xu,v)_\vep = -(u,Xv)_\vep\} =\{X\in \gl(V)|\, X^*=-X\}
    \\
    \fp(V) &:= \{X\in \gl(V)|\, (Xu,v)_\vep =  (u,Xv)_\vep\} = \{X\in \gl(V)|\, X^*=X\}.
    \end{aligned}
    \end{equation}

Let $\ft:=\ft(V)$ and $\fp:=\fp(V)$ for short.
The followings are immediate to check:
    \begin{equation}\label{eq: ftp}
    [\ft,\ft]\subset \ft, \ [\ft,\fp]\subset \fp,\ [\fp,\fp]\subset \ft.
    \end{equation}

Let $G(V):= \{g\in \GL(V)|\, (g v,g v')_\vep=(v,v')_\vep\ \mbox{for all $v,v'\in V$}\}$. Then $\ft=\gV$, where $\gV:=\Lie\, {G}(V)$. We have $\dim G(V)=\dim \ft=\frac12 k(k-\vep)$. So $\dim \fp=\frac12 k(k+\vep)$.

\begin{rk}
\label{rk: scalar}
If $V$ is a symplectic vector space of dimension 2 then $\fp(V)$ consists of scalars.
\end{rk}

If $\vep=-1$ then $G(V)$ is denoted by $\Sp(V)$ (called the symplectic group) with the Lie algebra $\fsp(V)$. If $\vep=+1$ then $G(V)$ is denoted by $\rO(V)$ (called the orthogonal group) with the Lie algebra $\fo(V)$. Unless no confusion arises, we shorten the notation as $\Sp,\fsp,\rO$ and $\fo$.


\subsection{}
\label{subsec: **}

It is direct to check that if $(\,,\,)_{V_1}=(\ ,\ )_{\vep}$ and $(\,,\,)_{V_2}=(\ ,\ )_{-\vep}$ for some $\vep\in \{-1,+1\}$ then $**=-\Id$. If both $(\,,\,)_{V_1}$ and $(\,,\,)_{V_2}$ are $(\ ,\ )_\vep$ for the same $\vep\in \{-1,+1\}$ then $**=\Id$. In particular if $V$ is symplectic and $W$ is orthogonal then for $i\in \Hom(W,V)$, $ii^*\in \fsp$ by \eqref{eq: t p} since $(ii^*)^*=i^{**}i^*=-ii^*$. Similarly $i^*i\in \fo$ for $i\in \Hom(W,V)$.


\subsection{Orbit-closedness and semisimplicity of quiver representations}
\label{subsec: the Orbit-closedness and the semisimplicity}
\begin{figure}[htbp]
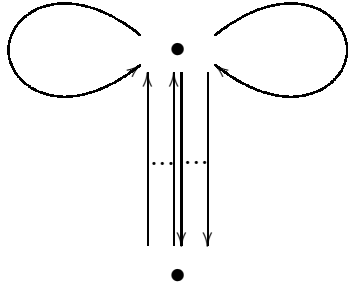

    \begin{equation*}
    \xy
    (0,40)*{\bullet}; (0,10)*{\bullet};
    %
    (-5,38)*{}; (-5,42)*{};  **\crv{(-15,50)&(-30,40)&(-15,30)} ?>*\dir{>};
    (-45,40)*{ };
    (5,38)*{}; (5,42)*{}; **\crv{(15,50)&(30,40)&(15,30)} ?>*\dir{>};
    (45,40)*{ };
    (-0.5,14)*{}; (-0.5,37)*{}  **@{-} ?>*\dir{>};
    (-4,14)*{}; (-4,37)*{} **@{-} ?>*\dir{>};
    (-2,25)*{...};
    (0.5,37)*{}; (0.5,14)*{}    **@{-} ?>*\dir{>};
    (4,37)*{}; (4,14)*{}  **@{-} ?>*\dir{>};
    (2.5,25)*{...};
    \endxy
    \end{equation*}
    \caption{a deframed quiver of the ADHM quiver}
    \label{fig: deframed quiver}
    \end{figure}

Let
    $$
    \bM:=\Hom(V,V)^{\oplus 2}\oplus \Hom(V,W)\oplus \Hom(V,W)
    $$
where $V$ and $W$ are vector spaces.
An element of $\bM$ is called an ADHM quiver representation. Let $W\cong \cc^N$ be any isomorphism. Then we have an obvious ($\GL(V)$-equivariant) linear isomorphism
    $$
    c\colon \bM\to \Hom(V,V)^{\oplus 2}\oplus\Hom(V,\cc)^{\oplus N}\oplus\Hom(V,\cc)^{\oplus N}.
    $$
We identify the target space of $c$ with the space of representations of the deframed quiver given in Fig.\ \ref{fig: deframed quiver} (cf. Crawley-Boevey's trick \cite[p.57]{CB2}). The number of arrows from the bottom vertex to the top is $N$, and the number of arrows from the top vertex to the bottom is also $N$.

Let $\vep\in \{-1,+1\}$. In the rest of paper we fix $V$ and $W$ as a $k$-dimensional vector space with $(\,,\,)_\vep$ and an $N$-dimensional vector space with $(\,,\,)_{-\vep}$ respectively.
Let
    $$
    \bN:=\left\{(B_1,B_2,i,j)\in \bM
    \middle\vert\, B_1=B_1^*,B_2=B_2^*,j=i^*\right\}.
    $$
It is clear that $\bN$ is isomorphic to the one we used in \S\ref{intro: subsec: ADHM} by the projection.

We use the notation $\stackrel{\perp}\oplus$ for a direct sum of vector spaces orthogonal with respect to the given nondegenerate forms.

It is well-known that $x\in \bM$ has a closed $\GL(V)$-orbit if and only if $c(x)$ is semisimple (cf.\ \cite{LP}). We have the corresponding result for $G(V)$ as below, whose proof is left as an exercise for the readers.

\begin{thm}\label{th: main 2}
Let $x:=(B_1,B_2,i,j)\in \bN$. Then the following properties are equivalent.
    \begin{itemize}
    \item[(a)] $c(x)$ is semisimple $($i.e. the direct sum of simple quiver representations$)$.
    \item[(b)] there exists a decomposition
        $$
        V=V^s\stackrel\perp\oplus \stackrel\perp\bigoplus_a V_a\stackrel\perp\oplus\stackrel\perp\bigoplus_b (V_b\oplus V_b')
        $$
    such that
    \begin{enumerate}
    \item $V_b$ and $V_b'$ are dual isotropic subspaces of $V$ for each index $b$;
    \item $c(B_1|_{V^s},B_1|_{V^s},i,j|_{V^s})$, $(B_1|_{V_a},B_2|_{V_a},0,0)$, $(B_1|_{V_b},B_2|_{V_b},0,0)$ and $(B_1|_{V_b'},B_2|_{V_b'},0,0)$ are simple quiver representations.
    \end{enumerate}
    \item[(c)] $G(V).x$ is closed in $\bN$.
    \end{itemize}
\end{thm}

We remark that $(B_1|_{V_b},B_2|_{V_b},0,0)$ and $(B_1|_{V_b'},B_2|_{V_b'},0,0)$ in the above statement are dual to each other.

Let us recall the two GIT stability conditions for the $\GL(V)$-actions on $\bM$ (cf.\ \cite[Chap.\ 2]{Lecture}).

\begin{defn}\label{defn: stable costable}
$(B_1,B_2,i,j)\in \bM$ is \textit{stable} (resp.\ \textit{costable}) if the following condition holds:
\begin{enumerate}
    \item
    (\textit{stability}) there is no subspace $S\subsetneq V$ such that $B_1(S)\subset S, B_2(S)\subset S$ and $\Image i\subset S$,
    \item
    (\textit{costability}) there is no nonzero subspace $T\subset V$ such that $B_1(T)\subset T, B_2(T)\subset T$ and $T\subset \Ker j$.
    \end{enumerate}
\end{defn}

Note that $x\in \bM$ is stable and costable if and only if $c(x)$ is simple when $W\neq0$.

Let $*_\bM\colon \bM\to \bM$, $(B_1,B_2,i,j)\mapsto (B_1^*,B_2^*,-j^*,i^*)$. Then $\bN$ is the $*_\bM$-invariant subspace of $\bM$. If $x\in \bM$ is stable (resp.\ costable) then $*_\bM(x)$ is costable (resp.\ stable). Therefore stability and costability are equivalent on $\bN$.
In particular, $x\in \bN$ is stable if and only if $c(x)$ is simple when $W\neq0$.


\subsection{ADHM description for framed $\Sp$-bundles and $\SO$-bundles}
\label{subsec: quiver}


We describe the framed $\Sp$-bundles and $\SO$-bundles in terms of ADHM data (cf.\ \cite{Do}, see also \cite{BrSan}).
Let $V:=\cc^k$ with $(\,,\,)_\vep$ and $W:=\cc^N$ with $(\,,\,)_{-\vep}$, where $\vep\in \{-1,+1\}$ and $N\ge1$. Let $\mu\colon \bN\to \fg(V)$ be the moment map given by $(B_1,B_2,i,j)\mapsto [B_1,B_2]+ij$.

Let $\mu\inv(0)^{\reg}$ be the (regular) locus of stable-costable quiver representations in $\mu\inv(0)$. We have $\mu\inv(0)^\reg\subset \mu\inv(0)^\sm$ since if $x\in \bN$ is stable then the differential $d\mu_x$ is surjective. By Theorem \ref{th: main 2} the image of $\mu\inv(0)^\reg$ of the GIT quotient map $\mu\inv(0)\to \mu\inv(0)\git G(V)$ is a Zariski open subset of $\mu\inv(0)\git G(V)$ and it is a $G(V)$-orbit space $\mu\inv(0)^\reg/G(V)$. In particular $\mu\inv(0)^\reg/G(V)$ is a smooth quasi-affine scheme. In fact it is also irreducible by \cite[Propositions 2.24 and 2.25]{BFG}.

\begin{rk}\label{rk: quiver description of Sp and So bundles}
By \S\ref{intro: subsec: ADHM2}, the scheme structures of $\cM_n^{\SO(N,\rr)}$ and $\cM_n^{\USp(N/2)}$ are defined as the quasi-affine scheme $\mu\inv(0)^\reg/G(V)$, where the relation between $n$ and $k$ are given there.\end{rk}

\begin{defn}
Let $\vep=-1$ (resp.\ $\vep=+1$). We say $x=(B_1,B_2,i,j)\in \mu\inv(0)$ is a \textit{$\SO$-datum} (resp.\ \textit{$\Sp$-datum}) if $[B_1,B_2]+ij=0$.
\end{defn}


\subsection{Stratification }
\label{subsec: stratification}

Let $\mu_k\colon \bN\to \fg(V)$, $(B_1,B_2,i,j)\mapsto [B_1,B_2]+ij$ (the holomorphic moment map).
To emphasize $k=\dim V$, we use $G_k$ for $G(V)$. It was denoted by $G_k'$ in \S\ref{sec: intro}.  Let $x:=(B_1,B_2,i,j)\in \mu_k\inv(0)$. Assume $G_k.x$ is closed in $\bN$. By Theorem \ref{th: main 2}, $(B_1|_{V^s},B_2|_{V^s},i,j|_{V^s})$ corresponds to a framed vector bundle, and $(B_1|_{V_a},B_2|_{V_a},0,0)$, $(B_1|_{V_b},B_2|_{V_b},0,0)$ and $(B_1|_{V_b'},B_2|_{V_b'},0,0)$ are all commuting pairs in $\fp(V)$. Let us focus on commuting pairs here. We simplify the situation: $B_1,B_2\in \fp(V)$ with $[B_1,B_2]=0$, where $V$ is orthogonal or symplectic. By the semisimplicity (Theorem \ref{th: main 2}), $B_1$ and $B_2$ are simultaneously diagonalizable. Therefore all $V_a,V_b$, and $V_b'$ of Theorem \ref{th: main 2} (b) are 1-dimensional. When $V$ is symplectic, $V_a$ does not appear and $B_1|_{V_b\oplus V_b'},B_2|_{V_b\oplus V_b'}$ are scalars, as $\fp(V)=\cc$ if $\dim V=2$. When $V$ is orthogonal, the index set of $b$ can be absorbed in the index set of $a$. We have the summary as follows.

\begin{thm}\label{th: main 4}
Let $S^n\aaa^2$ be the $n^{\mathrm{th}}$ symmetric product of $\aaa^2$.

$(1)$ Suppose $V$ is symplectic. Then there exists a canonical set-theoretical bijection
    $$
    \mu_k\inv(0)\git G_k=  \coprod_{0\le k'\le k}  \mu_{k'} \inv(0)^{\reg} / G_{k'} \times S^{\frac{k-k'}2}\aaa^2.
    $$

$(2)$ Suppose $V$ is orthogonal. Then there exists a canonical set-theoretical bijection
    $$
    \mu_k\inv(0)\git G_k=  \coprod_{0\le k'\le k}  \mu_{k'} \inv(0)^{\reg} / G_{k'} \times S^{k-k'}\aaa^2.
    $$
\end{thm}

Note that this stratification is nothing but the one of the Uhlenbeck space except the case $G=\SO(3)$ 
(cf.\ \cite[p.158, \S4.4.1]{DK}). For $G=\SO(3)$, $\mu_k\inv(0)\git G_k$ is different from the Uhlenbeck space, where the symmetric product is $S^{\frac{k-k'}4}\aaa^2$ by the formula \eqref{eq: rel of c2 and instanton number}.

Since $\mu_k\inv(0)^\reg/G_k$ is a free quotient (unless $\mu_k\inv(0)^\reg=\emptyset$), we have $\dim\mu_k\inv(0)^\reg/G_k=\dim \bN-2\dim G_k$. Using $\dim \fp-\dim\ft=k$ if $V$ is orthogonal, and $\dim \fp-\dim\ft=-k$ if $V$ is symplectic, we have
    \begin{equation}
    \label{eq: dim muk0}
    \dim\mu_k\inv(0)^\reg/G_k=\left\{\begin{array}{lll}
        k(N-2)& & \mbox{if $V$ is symplectic,}
        \\
        k(N+2)& & \mbox{if $V$ is orthogonal,}
        \end{array}
        \right.
    \end{equation}
whenever $\mu_k\inv(0)^\reg\neq\emptyset$.
Therefore we get the dimension of the strata. Hence by \S\ref{intro: subsec: ADHM}, we obtain the following.

\begin{lem}\label{lem: dim strata}
$(1)$ Assume $V$ is symplectic. Then $\mu_{k'} \inv(0)^{\reg} / G_{k'} \times S^{\frac{k-k'}2}\aaa^2$ is nonempty if and only if \textit{either} $N=3$ and $k'\in 4\zz_{\ge0}$, \textit{or} $N\ge4$ and $k'\in 2\zz_{\ge0}$. If it is  nonempty, it is of dimension $ k'(N-2)+(k-k')$

$(2)$ If $V$ is orthogonal then $\mu_{k'}\inv(0)^{\reg}/G_{k'} \times S^{k-k'}\aaa^2$ is nonempty and of dimension $ k'(N+2)+2(k-k')$. \qed
\end{lem}

\begin{rk}\label{rk: N3k4}
(1) Suppose $V$ is symplectic. By $\Spin(3)\cong \SU(2)$,  if $N=3$ and $k\in 4\zz_{\ge0}$ then $\mu_k\inv(0)^\reg/G_k$ is irreducible (\S\ref{intro: subsec: ADHM}). Hence, the strata of $\mu_4\inv(0)\git G_4$ are $\mu_4\inv(0)^\reg/G_4$ and $S^2\aaa^2$, both of which are irreducible varieties of dimension 4. See an alternative proof in Theorem \ref{th: N3k4} (3).

By $\Spin(4)\cong \SU(2)\times \SU(2)$, if $N=4$ and $k\in 2\zz_{\ge0}$ then $\mu_k\inv(0)^\reg/G_k$ has the $(k/2+1)$ irreducible components.

(2)
We will prove in Theorem \ref{th: main 5} that all the strata in Lemma \ref{lem: dim strata} (2) are irreducible normal.
\end{rk}


\section{Geometry of the moduli spaces of $\Sp$-data}
\label{sec: geometry of the moduli spaces of framed symplectic bundles}

Let $V$ be an orthogonal vector space of dimension $k$, and $W$ be a symplectic vector space of dimension $N\ge2$. Let $G:=G(V)=\rO(V)$.

The main theorem of this section is the following.

\begin{thm}
\label{th: main 5}
$\mu$ is normal, irreducible, flat, reduced and surjective. Hence $\mu\inv(0)$ is a normal irreducible variety of dimension $k(N+2)+k(k-1)/2$. And the Donaldson-Uhlenbeck space $\mu\inv(0)\git G$ is a normal irreducible variety of dimension $k(N+2)$, and $(\mu\inv(0)\git G)\setminus (\mu\inv(0)^\reg/G)$ is of codimension $\ge2$ in $\mu\inv(0)\git G$.
\end{thm}

The proof will appear in \S\ref{subsec: proof Theorem 5}.

\begin{cor}
The canonical embedding induces the equality of the ring of regular functions:
$\cc[\mu\inv(0)\git G]=\cc[\mu\inv(0)^\reg/G]$.
\end{cor}


\subsection{}\label{subsec: Panyushev}
Let $m\colon\fp\times \fp\to \ft$ by $(B_1,B_2)\mapsto [B_1,B_2]$.

\begin{thm}\label{th: Panyu}
For any $X\in \ft$, $m\inv (X)$ is an irreducible normal variety and a complete intersection in $\fp\times \fp$ of the $\dim\ft$ hypersurfaces. And the smooth locus of $m\inv(X)$ is the locus of $x\in m\inv(X)$ such that the differential $dm_x$ is surjective.
\end{thm}

\proof By \cite[(3.5)(1)]{Pa} and \cite[Theorem 3.2]{Pa}, $m\inv(X)$ is an irreducible complete intersection. See also \cite{Br} for a different proof. The normality is due to \cite[Corollary 4.4]{Pa}.

For $x\in m\inv(X)$, $x$ is a smooth point of $m\inv(X)$ if and only if $\dim T_xm\inv(X)=\dim m\inv(X)\, (=2\dim \fp-\dim \ft)$ if and only if $dm_x$ is surjective (cf. \cite[the proof of Proposition 4.2]{Pa}).\qed\vskip.3cm

For $B\in \gl$, $\gl^B=\{A\in \gl|\, [A,B]=0\}$. Let $\fp^B:= \gl^B\cap \fp$.

Let $\gl_l:= \{B\in \gl|\, \dim \gl^B=l\}$ and $\fp_l:= \{B\in \fp|\, \dim \fp^B=l\}$ ($l\in \zz_{\ge0}$).

\begin{lem}\label{lem: KR}
$($\cite[Prop.\ 5]{KR}$)$ For any $B\in \fp$, $\dim \fp^B-\dim \ft^B=k$. \qed \end{lem}

\begin{lem}\label{lem: local closed}
For any $l\in \zz_{\ge0}$,
$\fp_{\ge l}$ is a closed subvariety of $\fp$. And $\fp_k$ is Zariski open dense in $\fp$.
\end{lem}

\proof By \cite[p.7]{MFK}, for the conjugation action of $\GL(V)$ on $\gl$, the map $A\in\gl\mapsto \dim \GL(V)^A=\dim\gl^A$ is upper-semicontinuous.  Thus  $\gl_l$ is a locally closed subvariety of $\gl$. By \eqref{eq: ftp}, for any $B\in \fp$, $\gl^B=\fp^B\oplus \ft^B$. By Lemma \ref{lem: KR}, for $B\in \fp$, $B\in \fp_l$ if and only if $B\in \gl_{2l-k}$. Since $\fp$ is a closed subvariety of $ \gl$, $\fp_l=\fp\cap \gl _{2l-k}$ is a locally closed subvariety of $\fp$.

The second claim comes from the upper-semicontinuity, $\fp_k=\fp\cap \gl _{k}\neq \emptyset$ and $\gl_{l}=\emptyset$ for $l<k$. \qed\vskip.3cm

Let $\pi_i\colon\fp\times\fp\to \fp$ be the $i^{\mathrm{th}}$ projection $(i=1,2)$.

Note that for $X=0$, the singular locus of $m\inv(0)$ is of codimension $\ge3$ by \cite[p.6414, Lemma]{Br} (cf. \cite[Theorem (4.3)]{Pa}).

Let $M$ be a smooth variety, and $\phi\colon M\to \ft$ be a morphism.
The fibre-product $(\fp\times \fp)\times_\ft M$ is cut out by the $\dim\ft$-equations $[B_1,B_2]-\phi(x)=0$.

Let $S:=\{x\in\fp\times\fp|\, \mbox{$dm_x$ is surjective}\}$.

\begin{prop}\label{prop: main}
Let $M$ be a smooth irreducible variety and $\phi\colon M\to \ft$ be a morphism.
Then $(\fp\times \fp)\times_\ft M$ is an irreducible normal variety and a complete intersection in $\fp\times \fp\times  M$ of $\dim\ft$ hypersurfaces.
\end{prop}

\proof Let $\tm\colon(\fp\times \fp)\times_\ft M\to  M$ be the projection. Note that $\tm$ is surjective since so is $m$ by Theorem \ref{th: Panyu}. For each $x\in M$, $\tm\inv (x)\cong m\inv(\phi(x))$, which means any fibre dimension of $\tm$ is $2\dim\fp-\dim\ft$. Therefore $\dim(\fp\times \fp)\times_\ft M= 2\dim\fp-\dim\ft +\dim  M$ which proves $(\fp\times \fp)\times_\ft M$ is a complete intersection.

The complete intersection property implies that every irreducible component of $(\fp\times \fp)\times_\ft M$ has the same dimension.
Suppose $(\fp\times \fp)\times_\ft M$ is reducible.
There are two distinct irreducible components $Z,Z'$.
The restrictions to $Z,Z'$ of $\tm$ are dominant due to the equi-dimensionality of the fibers of $\tm$ and irreducibility of $M$.
For a generic element $x\in M$, $\tm\inv(x)\cap Z$ and $\tm\inv(x)\cap Z'$ are distinct closed subschemes in $\tm\inv(x)$ with the same dimension $2\dim\fp-\dim\ft$.
But this contradicts the fact that $m\inv(\phi(x))\, (\cong\tm\inv(x))$ is an irreducible scheme with dimension $2\dim\fp-\dim\ft$ (Theorem \ref{th: Panyu}).

By the upper-semicontinuity, $S$ is an open subvariety of $\fp\times \fp$. By Theorem \ref{th: Panyu}, $m\inv(X)\setminus S$ is of codimension $\ge2$ in $m\inv(X)$ for each $X\in \ft$. Thus $(\fp\times\fp)\setminus S$ is of codimension $\ge2$ in $\fp\times\fp$. By the smooth base change (\cite[III.10.1 (b)]{Hart}), $\tm_{S\times_\ft M}$ is a smooth morphism. By \cite[III.10.4]{Hart}, $S\times_\ft M$ is smooth. By Serre's criterion (\cite[Proposition II.8.23]{Hart}), we get the normality.
\qed\vskip.3cm

We have further description of the smooth locus of each fibre of $m$ as follows, which is not necessary in the proof of Theorem \ref{th: main 5}.

\begin{prop}\label{prop: smooth}
Let $X\in \ft$.

$(1)$
$m|\colon m\inv(X)\cap(\pi_1\inv(\fp_k)\cup \pi_2\inv(\fp_k)) \to \ft$ is a smooth morphism.

$(2)$ The codimension of
$m\inv(X)\setminus (\pi_1\inv(\fp_k)\cup \pi_2\inv(\fp_k))$ in $m\inv (X)$ is larger than 1.

\end{prop}

\proof (1) Let $L_B\colon\fp\to\ft$ by $A\mapsto [B,A]$, where $B\in \fp$. If $B\in\fp_k$ then $L_B$ is surjective since $\dim \ft+k=\dim\fp$. Since for $(B_1,B_2)\in \pi_1\inv(\fp_k)\cup \pi_2\inv(\fp_k)$, the differential $dm_{(B_1,B_2)}\colon\fp\times \fp\to \ft$ maps $(A_1,A_2)$ to $[B_1,A_2]+[A_1,B_2]$. Since $\ft=\Image L_{B_1}+ \Image L_{B_2}\subset \Image dm_{(B_1,B_2)}$, $dm_{(B_1,B_2)}$ is surjective. This proves (1).

(2) will be proven in Appendix \ref{app: B}.   \qed\vskip.3cm

\begin{rk}
Let $\tM:=(\fp\times \fp)\times_\ft M$, $\tS:=S\times_\ft M$ and $\tS':=(\pi_1\inv(\fp_k)\cup \pi_2\inv(\fp_k))\times_\ft M$. Then $\tS'\subset \tS\subset \tM$. By Proposition \ref{prop: main}, $\tS$ is a smooth variety such that $\tM\setminus\tS$ is of codimension $\ge2$ in $\tM$. In fact, we can strengthen the result further as follows: $\tS'$ is a smooth variety such that $\tM\setminus\tS'$ is of codimension $\ge2$.
By Proposition \ref{prop: smooth} (1) and the smooth base change, $\tS'$ is a smooth variety. By Proposition \ref{prop: smooth} (2),
$(\fp\times \fp) \setminus (\pi_1\inv(\fp_k)\cup \pi_2\inv(\fp_k))$ is of codimension $\ge2$ in $\fp\times \fp$. Thus $\tM\setminus \tS'$ is of codimension $\ge2$ in $\tM$.
\end{rk}


\subsection{The proof of Theorem \ref{th: main 5}}
\label{subsec: proof Theorem 5}

Let $M:= \Hom(W,V)$ and $\phi\colon M\to \ft$, $i\mapsto -ii^*$. By Proposition \ref{prop: main},
$\mu\inv(0)=(\fp\times \fp)\times_\ft M$ is an irreducible normal variety and a complete intersection in $\fp\times \fp\times  M$. By the method of associated cones \cite[II.4.2]{Kr}, $\mu$ is normal, irreducible, flat and reduced. Since $m$ is surjective, so is $\mu$.

By \cite[p.5]{MFK}, $\mu\inv(0)\git G$ is an irreducible normal variety. If one can show $\mu\inv(0)^\reg\neq\emptyset$ then by Lemma \ref{lem: dim strata} (2), $\mu\inv(0)^\reg/G$ is of the complement codimension $\ge2$, since $N\ge2$.

It remains to show $\mu\inv(0)^\reg\neq\emptyset$. Let $i\in \Hom(W,V)\setminus 0$. Let $X:=-ii^*$. Let $B_1\in \fp_k$. Let $a_n$ be the eigenvalues of $B_1$ ($n=1,2,...,k$). Let $V_{a_n}$ be the $a_n$-eigenspace. Let $p_n\colon V\to V_{a_n}$ be the projection. We may take $B_1$ so that $p_n(\Image(i))\neq0$ for all $n$. By Theorem \ref{th: Panyu}, for any $B_1\in \fp_k$, there exists $B_2\in \fp$ such that $[B_1,B_2]=X$. Now, $(B_1,B_2,i,i^*)\in \mu\inv(0)^\reg$, since $V=\sum_PP(B_1)i(W)$, where $P$ runs over one-variable polynomials. \qed\vskip.3cm


\section{Geometry of the moduli spaces of $\SO$-data: statements}
\label{sec: geometry of the moduli spaces of framed orthogonal bundles}

We assume that $V$ is symplectic and $W$ is orthogonal in the rest of the paper, except \S\ref{subsec: gen eigen} where preliminaries from linear algebra will be given. Let $k=\dim V$ and $N=\dim W$.

If $k=0$ then we have $\bN=0$.
If $N=1$ and $k\ge2$, $\mu\inv(0)^\reg=\emptyset$ as the right hand side of \eqref{eq: dim muk0} is negative. We will study the next simplest cases, $k=2$ (Theorem \ref{th: k2}), $N=2$ (Theorem \ref{th: N2}) and $(k,N)=(4,3)$ (Theorem \ref{th: N3k4}). $\mu$ is not flat nor irreducible in general unlike the case of the moduli spaces of $\Sp$-data.

In the following theorems, $\bP$ denotes the minimal nilpotent $\rO(W)$-orbit closure in $\fo(W)$. If $N=3$ or $N\ge5$ then $\bP$ is irreducible and normal. When $N=4$, $\bP$ has two irreducible components, which are isomorphic by the action of an element in $\rO(W)\setminus \SO(W)$. Each irreducible component is the closure of a $\SO(W)$-orbit. See \cite[Theorems 5.1.4 and 5.1.6]{CoMc}.

\begin{thm}\label{th: k2}
Let $k=2$. Then the following properties are verified.

$(1)$ $\mu\inv(0)$ is isomorphic to $\cc^2\times \{i\in \Hom(W,V)|\, ii^*=0\}$ by $(B_1,B_2,i,i^*)\mapsto (\tr B_1,\tr B_2,i)$.

$(2)$ If $N\ge 3$, $\mu$ is flat.

$(3)$ If $N\ge 5$, $\mu\inv(0)$ is irreducible and normal. Hence, $\mu$ is irreducible and normal. If $N=4$, $\mu\inv(0)$ is a reduced scheme and a union of two isomorphic irreducible components. If $N=3$, it is irreducible, but non-reduced.

$(4)$ If $N\le 3$, $\mu\inv(0)^\reg=\emptyset$. If $N\ge4$, $\mu\inv(0)^\reg=\mu\inv(0)^{\mathrm{sm}}$.

$(5)$ If $N\le 3$, the reduced scheme $\mu\inv(0)\git \Sp(V)_{\mathrm{red}}$ is isomorphic to $\cc^2$.

$(6)$ If $N\ge 4$, $\mu\inv(0)\git \Sp(V)$ is isomorphic to $\cc^2\times \bP$. Moreover the isomorphism restricts to $\mu\inv(0)^\reg/\Sp(V) \cong \cc^2\times (\bP\setminus0)$. \end{thm}

\begin{thm}\label{th: N2}
$(1)$
Let $N=2$ and $k\ge2$. Then $\mu\inv(0)$ is not a complete intersection. Hence $\mu$ is not flat.

$(2)$
Let $N=2$ and $k=4$. Then $\mu\inv(0)^\reg=\emptyset$.
\end{thm}

It is true that $\mu\inv(0)^\reg=\emptyset$ for $N=2$ and any $k\ge2$, as there is no nontrivial $\SO(2)$-instantons. But the author does not find its proof in terms of the ADHM description.

\begin{thm}\label{th: N3k4}
Let $N=3$ and $k=4$. Then we have the following assertions.

$(1)$ $\mu\inv(0)$ is a reduced complete intersection and a union of two irreducible components. Hence $\mu$ is flat.

$(2)$ One irreducible component of $\mu\inv(0)$ is the closure of $\mu\inv(0)^\reg$.

$(3)$ $\mu\inv(0)\git \Sp(V)$ is isomorphic to $\cc^2\times(\bP\sqcup_0 \bP)$. Moreover the isomorphism restricts to $\mu\inv(0)^\reg/ \Sp(V)\cong \cc^2\times(\bP\setminus 0)$.
\end{thm}

Here, $\bP\sqcup_0\bP:=(\bP\times\{0\})\cup (\{0\}\times \bP)$ in $\bP\times \bP$.
The proofs will appear in the subsequent sections.

Note that the first statement of Theorem \ref{th: k2} (4), Theorem \ref{th: N2} (2) and the second statement of Theorem \ref{th: N3k4} (3) are either obvious or well-known in the context of instantons on $S^4$ in \S\ref{intro: subsec: ADHM} (See also Remark \ref{rk: N3k4}).

The properties of the moment map $\mu$ in Theorem \ref{th: k2} (2), (3), Theorem \ref{th: N2} (1) and Theorem \ref{th: N3k4} (1) follow from the corresponding properties of $\mu\inv(0)$ by the method of associated cones \cite[II.4.2]{Kr}, as in the proof of Theorem \ref{th: main 5}.

The author is planning to study on flatness of $\mu$ for $N\ge4$ and normality of $\mu$ for $N\ge5$ in the near future. At least the following is true for small $k$ for a fixed $N$.

\begin{rk}\label{rk: Nge2k} It is known by \cite[Remark 11.3]{KP} that
$\rho\colon\Hom(W,V)\to \fsp(V)$, $i\mapsto ii^*$, is flat if $N\ge2k$. Moreover $\rho$ is normal and irreducible if $N\ge 2k+1$. By the base change argument used in the proof of Proposition \ref{prop: main}, the same is true for $\mu$. \end{rk}


\section{Kraft-Procesi's classification theory of nilpotent pairs}
\label{subsec: Kraft-Procesi's classification theory of nilpotent pairs}

In this section we review some geometry of $\Hom(W,V)$ following Kraft-Procesi \cite{KP}, which will be used in the proofs of the theorems from the previous section.


\subsection{Generalized eigenspaces and bilinear forms}\label{subsec: gen eigen}

Let $V$ be a vector space with $(\,,\,):=(\,,\,)_\vep$, where $\vep\in \{-1,+1\}$.

\begin{defn}
Let $X\in\gl(V)$. Let $V_a^n:=\Ker(X-a\Id)^n$ for $n\ge0$. We call $V_a:=\cup_{n\ge1}V_a^n$ the \textit{generalized $a$-eigenspace of} $X$.
\end{defn}

\begin{lem}
\label{lem: AA}
$(1)$ If $X\in \fp(V)$ and $a\neq b$ then $V_a\perp V_b$.

$(2)$ If $X\in \ft(V)$ and $a\neq -b$ then $V_a\perp V_b$.
\end{lem}

\proof (1) Let $m,n\ge0$. Let $v\in V_a^n$ and $w\in V_b^m$. We will show $(v,w)=0$. We use the induction on $m,n$. If $m$ or $n=0$, the assertion is obvious. Suppose the assertion is true for $m-1,n$ and $m,n-1$. Let $v':=(X-a)v\in V_a^{m-1}$ and $w':=(X-a)w\in V_b^{n-1}$. We have $(Xv,w)=a(v,w)$ as $(v',w)=0$ by our assumption. Similarly we also have $(v,Xw)=b(v,w)$. Thus we have $a(v,w)=b(v,w)$ which asserts $(v,w)=0$.

(2) The proof is same if one uses $(Xv,w)=-(v,Xw)$. \qed\vskip.3cm

Let $X\in \ft(V)$ or $\fp(V)$.
Let us define a bilinear form $|\,,\,|$ on $\Image (X)$ by $|Xv,Xv'|:=(v,Xv')$ (cf.\ \cite[\S4.1]{KP}).

\begin{lem}\label{lem: B}
$(1)$ If $X\in \ft(V)$ then $|\,,\,|$ on $\Image (X)$ is a nondegenerate bilinear form of type $-\vep$.

$(2)$ If $X\in \fp(V)$ then $|\,,\,|$ on $\Image (X)$ is a nondegenerate bilinear form of type $\vep$.
\end{lem}

\proof (1) is obvious, as noted in \cite[\S4.1]{KP}. (2) is also clear. \qed\vskip.3cm

\begin{prop}\label{prop: eigenspace dim ge 2}
Let $V$ be a $4$-dimensional symplectic vector space. Let $X\in\fp(V)$.
Then $X$ has an eigenspace of dimension $\ge2$.
\end{prop}

\proof By Lemma \ref{lem: AA}, \textit{either} $V=V_a\oplus V_b$ for some $a\neq b$, \textit{or} $V=V_a$ for some $a$. In the first case, $V_a$ and $V_b$ are 2-dimensional symplectic subspaces of $V$ and thus $X|_{V_a}=a$ and $X|_{V_b}=b$ by Remark \ref{rk: scalar}. In the second case, by Lemma \ref{lem: B}, $\Image(X-a)$ is a symplectic subspace of dimension 0 or 2. If  $\dim\Image(X-a)=0$ then $X=a$. If $\dim\Image(X-a)=2$ then   $X|_{\Image(X-a)}=a$ by Remark \ref{rk: scalar} as $X\in \fp(\Image(X-a))$. We are done. \qed\vskip.3cm

Let $W$ be a vector space with $(\,,\,)_{-\vep}$.
Recall that for $i\in \Hom(W,V)$, we have $ii^*\in \ft(V)$ and $i^*i\in \ft(W)$ (\S\ref{subsec: **}). Let $V_a$ (resp.\ $W_a$) be the generalized $a$-eigenspace of $ii^*$ (resp.\ $i^*i$).

\begin{lem}\label{lem: CC}
We have $ i^*(V_a) \subset W_a$ and $i(W_a)\subset V_a$.
Moreover if $a\neq0$ then $i$ and $i^*$ are isomorphisms between $ V_a$ and $W_a$.
\end{lem}

\proof For any $a\in\cc$ and $n\in \zz_{\ge0}$, we have $i^*(ii^*-a)^n= (i^*i-a)^ni^*$ and $(ii^*-a)^ni=i(i^*i-a)^n$. If $v\in V_a$ then $(i^*i-a)^ni^*v=i^*(ii^*-a)^nv=0$ for $n\gg0$. Similarly, if $w\in W_a$ then $(ii^*-a)^niw=i(i^*i-a)^nw=0$ for $n\gg0$. Therefore $i^*(V_a)\subset W_a$ and $i(W_a)\subset V_a$. This proves the first claim.

The restriction of $ii^*$ to $V_a$ is $a\Id_{V_a}$ plus a nilpotent endomorphism of $V_a$. So it is an isomorphism of $V_a$ if $a\neq0$. Similarly $i^*i$ gives an isomorphism of $W_a$. Therefore $i$ and $i^*$ give isomorphism between $V_a$ and $W_a$.
\qed\vskip.3cm


\subsection{Kraft-Procesi's results on nilpotent endomorphisms in $\ft$}
\label{subsec: Kraft-Procesi's classification of nilpotent pairs}

Let $\Sp:=\Sp(V),\mathrm{O}:=\mathrm{O}(W), \fsp:=\fsp(V)$ and $\fo:=\fo(W)$, for short. Let $X$ and $Y$ be nilpotent elements of $\fsp(V)$ and $\fo(V)$ respectively.

A nilpotent orbit $\Sp.X$ corresponds uniquely to a partition $\eta=(\eta_1,\eta_2,...)$  such that its transpose $(\heta_1,\heta_2,...)$ is given by $\heta_n:= \dim \Ker(X^n)/\Ker(X^{n-1})$. We denote $\eta$ by $\eta_X$. We represent $\eta_X$ by the Young diagram with the boxes replaced by $b$ (\textit{$b$-diagram}). For instance, a partition $(5,3,3,2,2,0,0,...)$ is represented by the $b$-diagram:
    $$
    \xy <1cm,0cm>:
    (1,1.2)*{bbbbb};
    (1,.8)*{bbb} ,(1,.4)*{bbb} , (1,0)*{bb} , (1,-0.4)*{bb}
    \endxy
    $$
For a nilpotent orbit $\rO.Y$, we use the symbol $a$ instead of $b$, and the Young diagram is called an \textit{$a$-diagram}.

\begin{rk}
\label{rk: SS}
By \cite[IV 2.15]{SS}, for $\eta_X=(\eta_1,\eta_2,...)$ (resp.\ $\eta_Y=(\eta_1,\eta_2,...)$),  $\#\{n|\, \eta_n=m\}$ is even for any odd $m$ (resp.\ for any even $m$).
\end{rk}

\begin{prop} $($\cite[Proposition 2.4]{KP}$)$
\label{prop: dim of Sp O}
For $\eta_X=\eta:=(\eta_1,\eta_2,...)$,
    $$
    \dim \Sp.X=\frac12\left(|\eta|^2+|\eta|-\sum\heta_n^2-\#\{n|\, \eta_n\, \mbox{is odd}\} \right) .
    $$
For $\eta_Y=\eta:=(\eta_1,\eta_2,...)$,
    $$
    \dim \mathrm{O}.Y= \frac12\left(|\eta|^2-|\eta|-\sum\heta_n^2+\#\{n|\, \eta_n\, \mbox{is odd}\}\right) .
    $$
\end{prop}

We say $\sigma\ge \eta$ for two partitions $\sigma=(\sigma_1,\sigma_2,...)$ and $\eta=(\eta_1,\eta_2,...)$ with $|\sigma|=|\eta|$ if $\sum_{1\le i\le j}\sigma_i\ge \sum_{1\le i\le j}\eta_i$ for any $j\ge1$ (\textit{dominance order} \cite[\S I.1]{Macdonald})

\begin{prop} $($\cite[Theorem 3.10]{He}$)$
\label{prop: Hesselink theorem}
Let $X$ and $X'$ be nilpotent elements in $\fsp$. Then $\eta_X\ge\eta_{X'}$ if and only if $\Sp.X'\subset \overline{\Sp.X}$. The same is true for $\fo$.
\end{prop}

\begin{cor}\label{cor: N}
$(1)$ If $\dim W\le 2$, there is no nonzero nilpotent orbit in $\fo$.

$(2)$ If $\dim W=3$, there exists a unique nonzero nilpotent orbit, and it is given by the $a$-diagram $aaa$.

$(3)$ If $\dim W\ge4$, the minimal nilpotent orbit is given by the $a$-diagram
    \begin{equation}
    \xy <1cm,0cm>:     (1, 0.6)*{ aa }, (1, 0.2)*{aa},   (1, -0.2)*{a }, (1, -0.6)*{\vdots}, (1, -1)*{a};
    \endxy
    \end{equation}
\end{cor}

\proof These assertions follow from Remark \ref{rk: SS} and Proposition \ref{prop: Hesselink theorem}.
\qed\vskip.3cm

Let us define two maps from $\Hom(W,V)$
    \begin{equation*}
    \begin{aligned}
    &
    \pi\colon\Hom(W,V)\to \fo,\ i\mapsto i^*i,
    \\
    &
    \rho\colon\Hom(W,V)\to \fsp, \ i\mapsto ii^*.
    \end{aligned}
    \end{equation*}

\begin{thm}\label{th: first fund th Inv th}
$($\cite[Theorem 1.2]{KP}$)$
$\pi$ and $\rho$ are the GIT quotient maps onto the images, by $\Sp$ and $\rO$ respectively.
\end{thm}

We review Kraft-Procesi's classification theory of nilpotent pairs in \cite{KP1} and \cite{KP}. A pair $(i,j)\in \Hom(W,V)\times\Hom(V,W)$ is a nilpotent pair if $ij$ is a nilpotent endomorphism. As in the case of nilpotent endomorphisms, a Young diagram plays an important role in the classification of the $\GL(V)\times \GL(W)$-orbits of nilpotent pairs.

\begin{defn}
\label{def: ab}
(\cite[\S\S4.2--4.3]{KP1})
By an $ab$-\textit{diagram}, we mean a Young diagram whose rows consists of alternating $a$ and $b$. E.g.,
    \begin{equation}
    \label{eq: ex of ab}
    \xy <1cm,0cm>:
    (1,1.2)*{ababa};
    (1,.8)*{aba} ,(1,.4)*{aba} , (1,0)*{ab} , (1,-0.4)*{ba}
    \endxy
    \end{equation}
\end{defn}

An $ab$-diagram $A$ gives a nilpotent pair as follows. Suppose the number of $a$ (resp.\ $b$) in $A$ is $k$ (resp.\ $N$).
Let us take any basis $\{b_1,b_2,...,b_k\}$ of $V$ (resp.\ $\{a_1,a_2,...,a_N\}$ of $W$). We replace all the $a$ and $b$ in $A$ by $a_i$ and $b_i$. We define a nilpotent pair $(i,j)$ such that $i$ maps $a_m$ to $b_n$ in the right adjacent position or 0 if there is no such $b_n$ and $j$ maps $b_m$ to $a_n$ or 0 similarly. In the above example of $ab$-diagram, we have a nilpotent pair as follows:
    \begin{equation}
    \xy <1cm,0cm>:
    (1,1.2)*{a_1\mapsto b_1\mapsto a_2\mapsto b_2\mapsto a_3\mapsto 0};
    (1,.8)*{a_4\mapsto b_3\mapsto a_5\mapsto 0} ,(1,.4)*{a_6\mapsto b_4\mapsto a_7\mapsto 0} , (1,0)*{a_8\mapsto b_5\mapsto 0} , (1,-0.4)*{b_6\mapsto a_9\mapsto 0}
    \endxy
    \end{equation}

In the above correspondence, an $ab$-diagram does not determine uniquely a nilpotent pair as the bases of $V$ and $W$ can be changed. Therefore up to the change of the bases, we have the bijective correspondence (\cite[\S4.3]{KP1})
    $$
    \begin{aligned}
    &
    \{ \mbox{$ab$-diagram with $\#a=\dim V$ and $\#b=\dim W$}\}
    \\ &
     \to
    \{(\GL(V)\times\GL(W)).(i,j)|\, \mbox{$(i,j)$ is a nilpotent pair}\}.
    \end{aligned}
    $$

Suppose $(i,i^*)$ is a nilpotent pair. From the $ab$-diagram of $(i,i^*)$, the $a$-diagram of $i^*i$ and the $b$-diagram of $ii^*$ are obtained by removing $b$ and $a$ respectively. For example,
    $$
    \xymatrix{
    {\xy <1cm,0cm>:
    (1,1.2)*{bb};
    (1,.8)*{b} ,(1,.4)*{b} , (1,0)*{b} , (1,-0.4)*{b}
    \endxy}
    &
    {\xy <1cm,0cm>:
    (1,1.2)*{ababa};
    (1,.8)*{aba} ,(1,.4)*{aba} , (1,0)*{ab} , (1,-0.4)*{ba}
    \endxy}
    \ar[r]^-{\pi} \ar[l]_-{\rho}
    &
    {\xy <1cm,0cm>:
    (1,1.2)*{aaa};
    (1,.8)*{aa} ,(1,.4)*{aa} , (1,0)*{a} , (1,-0.4)*{a}
    \endxy}
    }
    $$

For an $ab$-diagram we define
    \begin{equation}
    \Delta_{ab}:=\sum_{n:\, \mathrm{odd}} (\#\mbox{rows  of length $n$ starting with $a$})\cdot (\#\mbox{rows of length $n$ starting with $b$}).
    \end{equation}

\begin{thm}
\label{th: dim of Sp times O}
$($\cite[Theorem 6.5 and Proposition 7.1]{KP}$)$
Let $i\in \Hom(W,V)$.
\begin{enumerate}
    \item
    The orbits $(\Sp\times \mathrm{O}).i$ such that $X:=ii^*$ and $Y:=i^*i$ are nilpotent, are in 1-1 correspondence with the $ab$-diagrams whose rows are of one of the types $\alpha,\beta,\gamma, \delta,\ep$ from Table \ref{table: alpha beta} with $\# a=\dim W,\ \# b=\dim V$.
    \item
    For the $ab$-diagram associated to $i$,
        \begin{equation}
        \dim(\Sp\times \mathrm{O}).i= \frac12(\dim \Sp.X+\dim \OY+ \dim V.\dim W-\Delta_{ab}).
        \end{equation}
    \end{enumerate}
\end{thm}

\begin{table}[h]%
\caption{Rows of $ab$-diagrams}
\label{table: alpha beta}\centering %
\begin{tabular}{llllll}
Type & $\alpha_n$ & $\beta_n$ & $\gamma_n$ & $\delta_n$ & $\ep_n$   \vspace{0.1cm} \\ \hline  \vspace{0.1cm}
$ab$-diagram & $aba\cdots ba$   & $bab\cdots  ab$ &
    ${\xy <1cm,0cm>:
    (1,0.2)*{aba\cdots ba},
    (1,-0.2)*{aba\cdots ba}
    \endxy}$
&
    ${\xy <1cm,0cm>:
    (1,0.2)*{bab\cdots ab},
    (1,-0.2)*{bab\cdots ab}
    \endxy}$
&
    ${\xy <1cm,0cm>:
    (1,0.2)*{aba\cdots ab},
    (1,-0.2)*{bab\cdots ba}
    \endxy}$
    \vspace{0.1cm}
    \\ \hline \vspace{0.1cm}
$n$ &   &   & odd & even &  \vspace{0.1cm} \\ \hline \vspace{0.1cm}
$\# a$ & $2n+1$ & $2n-1$ & $2(n+1)$ & $2n$ & $2n$ \vspace{0.1cm} \\
$\# b$ & $2n$ & $2n$ & $2n$  & $2(n+1)$ & $2n$ \\ \hline \vspace{0.1cm}

\end{tabular}
\end{table}
%
%
%


\section{Moduli spaces of $\SO(N)$-data with $N\ge2$ and $k=2$}
\label{subsec: proof Theorem k2}

This section contains the proof of Theorem \ref{th: k2}.

Let $\dim V=k=2$ and $\dim W=N\ge2$.
Then $\fp:=\fp(V)$ consists of scalars by Remark \ref{rk: scalar}. Thus $(B_1,B_2,i,i^*)\in \mu\inv(0)$ implies $[B_1,B_2]=ii^*=0$, and $\mu\inv(0)\cong \cc^2\times  \rho\inv(0)$ by $(B_1,B_2,i,i^*)\mapsto (\tr(B_1),\tr(B_2),i)$. This proves (1).

Let $i\in \rho\inv(0)$.
Since $ii^*=0$, $b$ cannot appear twice in the same row in the $ab$-diagram of $i$. Looking at Table \ref{table: alpha beta}, we find that the $ab$-diagram of $i$ is one of the following:
    \begin{equation}\label{eq: ab when k2}
    \xy <1cm,0cm>:     (1, 0.6)*{ aba }, (1, 0.2)*{aba},   (1, -0.2)*{a }, (1, -0.6)*{\vdots}, (1, -1)*{a};
    \endxy
    \quad\quad\quad\quad\quad\quad
    \xy <1cm,0cm>:     (1, 0.6)*{ ab }, (1, 0.2)*{ba},   (1, -0.2)*{a }, (1, -0.6)*{\vdots}, (1, -1)*{a};
    \endxy
    \quad\quad\quad\quad\quad\quad
    \xy <1cm,0cm>:     (1, 0.6)*{ b }, (1, 0.2)*{b},   (1, -0.2)*{a }, (1, -0.6)*{\vdots}, (1, -1)*{a};
    \endxy
    \end{equation}
where the left-most one actually occurs only when $N\ge4$. We denote elements of $\Hom(W,V)$ corresponding to the above $ab$-diagrams by $i_1,i_2$ and $i_3(=0)$ respectively, where $i_1$ does not exist unless $N\ge4$. By Theorem \ref{th: dim of Sp times O}, $\dim (\Sp\times \rO).i_1=2N-3$, $\dim (\Sp\times \rO).i_2=N$ and $\dim (\Sp\times \rO).i_3=0$. Therefore $\rho\inv(0)$ is of the expected dimension. By the method of associated cones \cite[II.4.2]{Kr}, $\rho$ is flat. Thus $\mu$ is flat. This proves (2). See Remark \ref{rk: Nge2k} for the flatness when $N\ge4$.

Suppose $N=3$. Then $\rO=\SO\sqcup -\SO$. Since $-\Id_V\in \Sp$ and $-\Id_W.i =-\Id_V.i $ for any $i\in \Hom(W,V)$, $(\Sp\times \rO).i_2 \,(=(\Sp\times \SO).i_2)$ is irreducible and its Zariski closure is $\rho\inv(0)$. Let us check $\rho\inv(0)$ is not reduced. By \cite[Remark 11.4]{KP}, $\rho\inv(0)^\sm=\{i\in \rho\inv(0)|\, \mbox{$i$ is surjective}\}$. By \eqref{eq: ab when k2} any $i\in \rho\inv(0)$ is not surjective, which implies $\rho\inv(0)$ is not reduced.

If $N=4$ then by \cite[\S11.3]{KP}, $\rho\inv(0)$ is reduced and consists of two irreducible components which are isomorphic by the action of an element of $\rO\setminus\SO$.

If $N\ge5$, \cite[\S11.3]{KP} asserts $\rho\inv(0)$ is normal and irreducible.
This proves (3).

Let $x:=(B_1,B_2,i,i^*)\in\mu\inv(0)$.
By the $ab$-diagrams in
\eqref{eq: ab when k2}, $\Ker(i_2^*)$ and $\Ker(i_3^*)$ are nonzero. Any nonzero vector in $\Ker(i_2^*)$ or $\Ker(i_3^*)$ is a common eigenvector of scalars $B_1$ and $B_2$. Thus, if $i=i_2$ or $i_3$ then $x$ is not costable. Since $\Ker(i_1^*)=0$, if $i=i_1$ then $x$ is costable. So $\mu\inv(0)^\reg=\emptyset$ if $N\le 3$. If $N\ge4$,  $\mu\inv(0)^\reg=(\Sp\times\rO).i_1$.

We describe the smooth locus of $\mu\inv(0)=\cc^2\times \rho\inv(0)$ when $N\ge4$. Since $T_x\rho\inv(0)=\Ker d\rho_x$ and $\dim \rho\inv(0)=\dim \Hom(W,V)-\dim \Sp$, $\rho\inv(0)^\sm$ is the locus of $x$ such that $d\rho_x$ is surjective. Thus it is the locus of $x$ such that $\rho$ is smooth at $x$. By \cite[Remark 11.4]{KP}, $\rho\inv(0)^\sm=\{i\in \rho\inv(0)|\, \mbox{$i$ is surjective}\}$. Among $i_1,i_2$ and $i_3$, only $i_1$ is surjective. This proves (4).

Since $\Sp$ acts trivially on $\fp$, we have $\mu\inv(0)\git \Sp\cong \cc^2\times (\rho\inv(0)\git \Sp)$.
By Theorem \ref{th: first fund th Inv th} and the $ab$-diagrams in \eqref{eq: ab when k2}, the reduced scheme
    $$
    \rho\inv(0)\git \Sp{}_{\mathrm{red}}\cong\left\{\begin{array}{lll}
        0 & \mbox{if $N\le3$}
        \\
        \bP & \mbox{if $N\ge4$.}
        \end{array}
        \right.
    $$
This proves the statement on $\mu\inv(0)\git\Sp$ in (5) and (6). By Theorem \ref{th: first fund th Inv th}  and the $a$-diagrams coming from the $ab$-diagrams in \eqref{eq: ab when k2}, we have $\overline{(\Sp\times \rO).i_2}\git\Sp =0$. This proves the second claim of (6).
\qed\vskip.3cm


\section{Moduli spaces of $\SO(2)$-data }
\label{subsec: proof Theorem N2}

This section contains the proof of Theorem \ref{th: N2} (1).
The proof of Theorem \ref{th: N2} (2) will appear in \S\ref{sec: proof of Theorem N2 (2)}.

Let $\dim V=k\ge2$ and $\dim W=N=2$.
Let $m\colon \fp\times\fp\to \ft$, $(B,B')\mapsto [B,B']$.

We will show there exists a subvariety $X$ in $\mu\inv(0)$ of dimension $>$ $2\dim \fp -\dim \ft + \dim \Hom(W,V)$. 
Let $ X:=\{(B,B',i,i^*)|\, [B,B']=0=ii^*\}=m\inv(0)\times \rho\inv(0) \subset \mu\inv(0)$.

Let us estimate $\dim m\inv(0)$. Let $\fp^{(e)}:=\{B\in\fp|\, \mbox{$B$ has distinct $e$ eigenvalues}\}.$ As shown in Appendix \ref{app: B}, $\fp^{(\le e)}$ is a closed subvariety of $\fp$. By Lemma \ref{lem: AA} (1), if $e>k/2$, $\fp^{(e)}=\emptyset$. With respect to a symplectic basis of $V$, $\diag(a_1,a_2,...,a_{k/2})^{\oplus2}\in  \fp^{(k/2)}$, where $a_1,a_2,...,a_{k/2}$ are all distinct. Therefore $\fp^{(k/2)}$ is a Zariski dense open subset of $\fp$. Let $p\colon m\inv(0)\to \fp$ be the first projection. It is clear that $p$ is surjective.  And $p\inv(B_0)\cong \fp^{B_0}$ for any $B_0\in \fp$ (see \S\ref{subsec: Panyushev} for the notation). If $B_0\in \fp^{(k/2)}$ then we claim $\dim\fp^{B_0}=k/2$. We have the (2-dimensional) eigenspace decomposition $V=\bigoplus_{i=1}^{k/2}V_{a_i}$ of $B_0$. Since $\gl^{B_0}= \bigoplus_{i=1}^{k/2}\gl(V_{a_i})$ we have $\fp^{B_0}=\bigoplus_{i=1}^{k/2}\fp(V_{a_i})=\bigoplus_{i=1}^{k/2}\cc$.

By the claim we obtain an estimate:
    $$
    \dim m\inv(0)\ge \frac k2+\dim \fp.
    $$

Let us compute $\dim\rho\inv(0)$. Let $i\in \rho\inv(0)$. Since $\dim W=2$ and $i^*i$ is nilpotent, $i^*i=0$ (Corollary \ref{cor: N} (1)). By an argument similar to the one we used in \S\ref{subsec: proof Theorem k2}, the  $ab$-diagram of any nonzero $i$ is
    $$
    \xy <1cm,0cm>:     (1, 0.6)*{ ab }, (1, 0.2)*{ba},   (1, -0.2)*{b }, (1, -0.6)*{\vdots}, (1, -1)*{b};
    \endxy
    $$
Thus $\rho\inv(0)$ is the Zariski closure of $(\Sp\times \mathrm{O}).i$ and $\dim\rho\inv(0)= k$ by Theorem \ref{th: dim of Sp times O} (2).

To sum up, $\dim X\ge\frac32 k+\dim \fp >2\dim \fp-\dim \ft+2k$ because $\dim \ft-\dim\fp=k$, which means $\mu\inv(0) $ is not a complete intersection in $\bN$.
\qed


\section{Moduli spaces of $\SO(3)$-data with $k=4$}
\label{sec: the moduli spaces of framed orthogonal bundles of rank 3 and k4}
This section contains the proof of Theorem \ref{th: N3k4}.

Let $\dim V=k=4$ and $\dim W=N=3$.
Let $\SO:=\SO(W) $ for short.


\subsection{Description of $\mu\inv(0)$}
\label{subsec: description of mu inv 0 N3k4}

Let $\fp':=\{D\in \fp|\, \tr D=0\}$. Let $\tX:=\{(B_1,B_2,i,i^*)\in \mu\inv(0)|\, B_1,B_2\in \fp'\}$. Then
    \begin{equation*}
    \begin{aligned}
    &
    \mu\inv(0)\cong \cc^2\times \tX ,
    \\
    &
    (B_1,B_2,i,i^*)\mapsto \left((\tr(B_1),\tr(B_2)),(B_1-\frac14\tr(B_1)\Id_V,B_2-\frac14\tr(B_2)\Id_V,i,i^*) \right).
    \end{aligned}
    \end{equation*}
Let $\bN':=\fp'{}^{\oplus2}\oplus \{(i,i^*)|\, i\in \Hom(W,V)\}$. Then
$\tX$ is defined by $\mu=0$ in $\bN'$. Therefore it is enough to show the corresponding statements for $\tX$ to prove Theorem \ref{th: N3k4}.

Let $\tX_1:=\{(B_1,B_2,i,i^*)\in \tX|\, [B_1,B_2]=0\}$. Then $\tX_1=\{(B_1,B_2)\in \fp'\times \fp'|\, [B_1,B_2]=0\}\times \rho\inv(0).$

\begin{lem}\label{lem: 121129 Claim (1)}
$\rho\inv(0)$ consists of two irreducible strata which are $(\Sp\times \rO)$-orbits of dimension 5 and 0 respectively.
Hence $\rho\inv(0)$, as a cone, is irreducible.
\end{lem}

\proof By a similar argument as in \S\ref{subsec: proof Theorem k2}, the possible $ab$-diagrams of $i$ of $\rho\inv(0)$ are
    \begin{equation}
    \label{eq: ab IV}
    {\xy <1cm,0cm>:
    (1,0.8)*{ab };
    (1,0.4)*{ba };
    (1,0 )*{b  };
    (1,-0.4)*{ b };
    (1,-0.8)*{a  };
    \endxy} \quad\quad\quad\quad\quad\quad
    {\xy <1cm,0cm>:
    (1,1.2)*{b };
    (1,0.8)*{b  };
    (1, 0.4)*{  b };
    (1, 0)*{  b };
    (1, -0.4)*{  a };
    (1, -0.8)*{  a };
    (1, -1.2)*{  a };
    \endxy}
    \end{equation}
By Theorem \ref{th: dim of Sp times O}, $\rho\inv(0)$ consists of two $(\Sp\times \rO)$-orbits associated to the above $ab$-diagrams of dimension 5 and 0 respectively.

Since $-\Id_W\in \rO\setminus \SO$ and $-\Id_V\in \Sp$, $(\Sp\times \rO).i=(\Sp\times \SO).i$ for any $i\in \Hom(W,V)$. So we have irreducibility.
\qed\vskip.3cm

Let $e_1,e_2,e_3,e_4$ be a basis of $V$ such that $(e_1,e_2)_V=(e_3,e_4)_V=1\ \mbox{and}\ (e_l,e_m)_V=0$ for other $l,m$ with $l\le m$.
Let $f_1,f_2,f_3$ be an orthogonal basis of $W$ so that $(f_i,f_j)=\delta_{ij} .$

Let $I$ be the $2\times 2$ identity matrix.
Let
    \begin{equation}
    J:=\begin{pmatrix} 0 &-1
    \\ 1 & 0
    \end{pmatrix}
    ,\
    H:=\begin{pmatrix} 1 &0
    \\ 0 & -1
    \end{pmatrix},\
    X:=\begin{pmatrix} 0 & 1
    \\ 0 & 0
    \end{pmatrix}
    \ \mbox{and}\
    Y:=\begin{pmatrix} 0 & 0
    \\ 1 & 0
    \end{pmatrix}
    .
    \end{equation}

We identify $\gl(V)=\Mat_4$ with respect to $e_1,...,e_4$. Then we can write the elements of $\ft$ and $\fp$ in matrix forms in $2\times 2$ subminors:
    \begin{equation}
    \label{eq: ft fp}
    \begin{aligned}
    \ft&:=\left\{
        \begin{pmatrix}
        P & JQ
        \\
        JQ^t  & S
        \end{pmatrix}
        \in
        \gl(V)\middle|\,
        \tr(P)=\tr(S)=0   \right\},
        \\
    \fp&:=\left\{
        \begin{pmatrix}
        aI & JR
        \\
        -JR^t  & bI
        \end{pmatrix}
        \in
        \gl(V) \middle| \, a,b\in\cc \right\}.
    \end{aligned}
    \end{equation}
where $P,Q,R,S\in \Mat_2$.
Note that
    $$\begin{pmatrix}
        0 & H
        \\
        -H & 0
        \end{pmatrix},\
        \begin{pmatrix}
        0 & X
        \\
        -X & 0
        \end{pmatrix} ,\
        \begin{pmatrix}
        0 & I
        \\
        I & 0
        \end{pmatrix} \
        \mbox{and}\
        \begin{pmatrix}
        0 & Y
        \\
        -Y & 0
        \end{pmatrix}
        $$
are elements of $\fp'$ by letting $R=-JH$, $R=-JX$, $R=-JY$ and  $R=-J$ respectively. Therefore we obtain a basis of $\fp'$ as
    \begin{equation}
    \label{eq: basis decomposition of fp}
    v_1:= \frac12
        \begin{pmatrix}
        I & 0
        \\
        0 & -I
        \end{pmatrix},
    v_2:= \begin{pmatrix}
        0 & H
        \\
        -H & 0
        \end{pmatrix},
    v_3:=
        \begin{pmatrix}
        0 & X
        \\
        -X & 0
        \end{pmatrix} ,
    v_4:=
        \begin{pmatrix}
        0 & Y
        \\
        -Y & 0
        \end{pmatrix} ,
    v_5:= \frac12
        \begin{pmatrix}
        0 & I
        \\
        I & 0
        \end{pmatrix}.
    \end{equation}
By direct computation, the Lie brackets of pairs of basis elements of $\fp'$ are
    \begin{equation}\label{eq: list Lie bracket}
    \begin{aligned}
    &
    [v_1,v_2]= \begin{pmatrix} 0 & H \\ H & 0 \end{pmatrix},\
    [v_1,v_3]= \begin{pmatrix} 0 & X \\ X & 0 \end{pmatrix},\
    [v_1,v_4]= \begin{pmatrix} 0 & Y \\ Y & 0 \end{pmatrix},\
    [v_1,v_5]= \frac12\begin{pmatrix} 0 & I \\ -I & 0 \end{pmatrix},\
    \\
    &
    [v_2,v_3]= -2\begin{pmatrix} X & 0 \\ 0 & X \end{pmatrix},\
    [v_2,v_4]= 2 \begin{pmatrix} Y & 0 \\ 0 & Y \end{pmatrix},\
    [v_2,v_5]= \begin{pmatrix} H & 0 \\ 0 & -H \end{pmatrix},\
    \\
    &
    [v_3,v_4]= -\begin{pmatrix} H & 0 \\ 0 & H  \end{pmatrix},\
    [v_3,v_5]= \begin{pmatrix} X & 0 \\ 0 & -X \end{pmatrix},\
    [v_4,v_5]= \begin{pmatrix} Y & 0 \\ 0 & -Y  \end{pmatrix}.
    \end{aligned}
    \end{equation}
These form a basis of $\ft$.

Let us define a linear map $F\colon \wedge^2\fp'\to \ft$ by setting $F(v_i\wedge v_j):=[v_i,v_j]$, where $1\le i<j\le 5$. Then we have a commuting diagram
    \begin{equation}
    \label{eq: ft wedge fp}
    \xymatrix{
        \fp'\times \fp' \ar[r]^-\omega \ar[rd]_{[\ ,\ ]} & \wedge^2\fp' \ar[d]^F
        \\
        &   \ft.
        }
    \end{equation}
By \eqref{eq: list Lie bracket},  $\Image(F)$ contains the basis of $\ft$ and $\dim \ft=\dim \wedge^2\fp'=10$. Thus $F$ is an isomorphism. In particular, for $B_1,B_2\in \fp'$, $[B_1,B_2]=0$ if and only if $B_1\wedge B_2=0$. This proves the following lemma.

\begin{lem}
\label{lem: 121129 Claim (2)}
$\{(B_1,B_2)\in \fp'\times \fp'|\, [B_1,B_2]=0\}= \{(aB,bB)|\, a,b\in \cc,\ B\in \fp'\}$ and it is irreducible. \qed
\end{lem}

\begin{cor}
\label{cor: tX1}
$(1)$ $\tX_1$ is irreducible of dimension 11.

$(2)$ Any element of $\tX_1$ is an unstable quiver representation.
\end{cor}

\proof
(1) follows from Lemmas \ref{lem: 121129 Claim (1)} and \ref{lem: 121129 Claim (2)}.

(2) The non-costability amounts to the existence of a common eigenvector $v$ of $B_1$ and $B_2$ such  that $i^*(v)=0$. From \eqref{eq: ab IV}, we have $\dim \Ker (i^*)\ge3$. Let $B_1=a_1B$ and $B_2=a_2B$ for some $a_1,a_2\in \cc$ and $B\in\fp'$ (Lemma \ref{lem: 121129 Claim (2)}). By Proposition \ref{prop: eigenspace dim ge 2}, we have an eigenspace of $B$ of dimension $\ge2$. By the dimension reason it has a nonzero vector contained in the eigenspace and $\Ker(i^*)$. \qed\vskip.3cm

On the other hand, if $\sigma\in \Image\,\omega$ is nonzero then we have
    \begin{equation}
    \label{eq: Sl2}
    \omega\inv(\sigma)\cong \SL(2).
    \end{equation}
since $\Image\,\omega\setminus 0$ is the set of 2-dimensional subspaces $S$ of $\fp'$ with a volume form of $S$.
This isomorphism can be described in  more detail as follows.  Define an $\SL(2)$-action on $\fp'\times \fp'$ as
    \begin{equation}
    \label{eq: Sl2 V V}
    \begin{pmatrix} a & b \\ c& d \end{pmatrix}.(B_1,B_2):=(aB_1+bB_2,cB_1+d B_2).
    \end{equation}
The $\SL(2)$-action on ${\omega\inv(\wedge^2\fp'\setminus0)}$ is free.
The $\SL(2)$-action on $\wedge^2\fp'$ is trivial and $\omega|_{\omega\inv(\wedge^2\fp'\setminus0)}$ is $\SL(2)$-equivariant. Now \eqref{eq: Sl2} is nothing but the identification of a free $\SL(2)$-orbit.

\begin{defn}\label{def: G}
Let $G:=\SL(2)\times \Sp\times \rO$. Then we have a $G$-action on $\tX$ by
    $$
    (g_1,g_2,g_3).((B_1,B_2),i,i^*):= (g_1.(g_2B_1g_{2}\inv,g_2B_2g_{2}\inv),g_2ig_3\inv,(g_2ig_3\inv)^*).
    $$
\end{defn}

\begin{lem}\label{lem: 121206 claim pre}
Let $B_1,B_2\in \fp'$. Then $([B_1,B_2])^2$ is a scalar endomorphism.
\end{lem}

\proof This follows from a tedious, but direct computation.
\qed\vskip.3cm

Let $\tX_2:=\tX\setminus \tX_1=\{(B_1,B_2,i,i^*)\in \tX |\, [B_1,B_2]\neq0\}$. Let $p\colon \tX_2\to \Hom(W,V)$ be the projection.

\begin{cor}
\label{cor: 121206 claim}
Let $i\in p(\tX_2)$. Then $(ii^*)^2=0$.
\end{cor}

\proof  Write $i=p(B_1,B_2,i,i^*)$ for some $B_1,B_2\in \fp'$.
By Lemma \ref{lem: 121206 claim pre}, $[B_1,B_2]^2=(ii^*)^2$ is a scalar endomorphism. Since $\rank \, i\le 3$, the scalar is 0. \qed

\begin{cor}
\label{cor: 121206 cor}
Let $i\in p(\tX_2)$. Then the $ab$-diagram of $i$ is one of the following:
    \begin{equation}
    \label{eq: ab III}
    (\I)\quad
    {\xy <1cm,0cm>:
    (1,0.4)*{bab };
    (1,0)*{bab };
    (1,-0.4)*{ a };
    \endxy} \quad\quad\quad
    (\II)\quad
    {\xy <1cm,0cm>:
    (1,0.4)*{ababa };
    (1,0)*{b };
    (1,-0.4)*{ b };
    \endxy} \quad\quad\quad
    (\III)\quad
    {\xy <1cm,0cm>:
    (1,0.4)*{bab };
    (1,0)*{ba };
    (1,-0.4)*{ ab };
    \endxy} \quad\quad\quad
    (\IV)\quad
    {\xy <1cm,0cm>:
    (1,0.8)*{bab };
    (1,0.4)*{b };
    (1, 0 )*{ b };
    (1, -0.4 )*{ a };
    (1, -0.8 )*{ a };
    \endxy}
    \end{equation}
\end{cor}

\proof
Since $(ii^*)^2=0$ by Corollary \ref{cor: 121206 claim} and $ii^*\neq0$, the maximal length of rows of the $b$-diagram of $ii^*$ is 2. Thus the $b$-diagram of $ii^*$ is one of the following:
    \begin{equation}
    \xy <1cm,0cm>:
    (1,.2)*{bb};
    (1,-0.2)*{bb}
    \endxy \quad\quad\quad
    \xy <1cm,0cm>:
    (1,.4)*{bb};
    (1,-0 )*{b};
    (1,-0.4 )*{b};
    \endxy
    \end{equation}
From Table \ref{table: alpha beta}, we get the list of $ab$-diagrams as above. \qed\vskip.3cm

\begin{lem}\label{lem: 121206 Claim (1)(2)}
$(1)$ Let $B_1:=v_5$ and $B_2:=v_1-\frac12 v_2$. Then $[B_1,B_2]$ is a nilpotent matrix whose $b$-diagram is
    \begin{equation}
    \xy <1cm,0cm>:
    (1,.2)*{bb};
    (1,-0.2)*{bb}
    \endxy
     \end{equation}
$(2)$ Let $B_1:=v_3$ and $B_2:=\sqrt{-1}v_1- v_5$. Then $[B_1,B_2]$ is a nilpotent matrix whose $b$-diagram is
    \begin{equation}
    \xy <1cm,0cm>:
    (1,.4)*{bb};
    (1,-0 )*{b};
    (1,-0.4 )*{b};
    \endxy
    \end{equation}
\end{lem}

\proof (1) By direct calculation we have
    $$
    [B_1,B_2]=\frac12\begin{pmatrix}
        1 & 0 & -1 & 0 \\ 0 & -1 & 0 & -1 \\
        1 & 0 & -1 & 0 \\ 0 & 1 & 0 & 1
        \end{pmatrix}
    $$
which maps
    \begin{equation}\label{eq: i3i3*}
    e_1\mapsto \frac12(e_1+e_3),\ \frac12(e_1+e_3)\mapsto0,\ e_2\mapsto -\frac12(e_2-e_4),\  \frac12(e_2-e_4)\mapsto0 .
    \end{equation}

(2) By direct calculation we have
    $$
    [B_1,B_2]= \begin{pmatrix}
        0 & -1 & 0 & -\sqrt{-1} \\ 0 & 0 & 0 & 0 \\
        0 & -\sqrt{-1} & 0 & 1 \\ 0 & 0 & 0 & 0
        \end{pmatrix}
    $$
which maps
    \begin{equation}\label{eq: i4i4*}
    e_2\mapsto -(e_1+\sqrt{-1}e_3),\ e_1+\sqrt{-1}e_3\mapsto0,\ e_1\mapsto 0,\ e_2+\sqrt{-1}e_4\to 0 .
    \end{equation}
\qed\vskip.3cm

We denote $B_1,B_2$ in Lemma \ref{lem: 121206 Claim (1)(2)} (1) by $B_1^{(\I)},B_2^{(\I)}$ respectively.
We denote $B_1,B_2$ in Lemma \ref{lem: 121206 Claim (1)(2)} (2) by $B_1^{(\II)},B_2^{(\II)}$ respectively. Let $B_n^{(\III)}$ and $B_n^{(\IV)}$ be $B_n^{(\II)}$ ($n=1,2$).
For each $A=\I,\II,\III,\IV$ there exists $i^{(A)}\in \Hom(W,V)$ whose $ab$-diagram is of type $(A)$ in \eqref{eq: ab III} and $[B_1^{(A)},B_2^{(A)}]+{i^{(A)}}{i^{(A)}}^*=0$ by Theorem \ref{th: dim of Sp times O} (1).

Let us define $i^{(\I)}$ in a matrix form. We take a basis $\{\te_1,\te_2,\te_3,\te_4\}$ of $V$ by setting $\te_1:=e_1$, $\te_2:=e_2-e_4$, $\te_3:=e_1+e_3$ and $\te_4:=e_4$. We take a basis $\{\tf_1,\tf_2,\tf_3\}$ of $W$ by setting $\tf_1:=\frac1{\sqrt{2}}(f_1+\sqrt{-1}f_2)$, $\tf_2:=\frac1{\sqrt{2}}(f_1-\sqrt{-1}f_2)$ and $\tf_3:=f_3$. Let
    $$
    i^{(\I)}:=\frac1{\sqrt{-2}}\begin{pmatrix}
        0 & 0 & 0 \\
        0 & -1 & 0 \\
        1 & 0 & 0 \\
        0 & 0 & 0
        \end{pmatrix}
    $$
with respect to $\{\te_1,\te_2,\te_3,\te_4\}$ and $\{\tf_1,\tf_2,\tf_3\}$.

Let us define $i^{(\II)}$, $i^{(\III)}$ and $i^{(\IV)}$ in matrix forms respectively. We take a basis $\{\te_1',\te_2',\te_3',\te_4'\}$ of $V$ by setting $\te_1':=e_2$, $\te_2':=-(e_1+\sqrt{-1}e_3)$, $\te_3':=e_3$ and $\te_4':=-\sqrt{-1}e_2+e_4$. We take a basis $\{\tf_1',\tf_2',\tf_3'\}$ of $W$ by setting $\tf_1':=\frac1{\sqrt{2}}(f_1+\sqrt{-1}f_3)$, $\tf_2':=f_2$ and $\tf_3':=\frac1{\sqrt{2}}(f_1-\sqrt{-1}f_3)$. Let
    $$
    i^{(\II)}:= \begin{pmatrix}
        1 & 0 & 0 \\
        0 & -1 & 0 \\
        0 & 0 & 0 \\
        0 & 0 & 0
        \end{pmatrix},
    \quad
    i^{(\III)}:=\begin{pmatrix}
        0 & 0 & 0 \\
        0 & -1 & 0 \\
        1 & 0 & 0 \\
        0 & 0 & 0
        \end{pmatrix}
    \quad\mbox{and}\quad
    i^{(\IV)}:= \begin{pmatrix}
        0 & 0 & 0 \\
        0 & -1 & 0 \\
        0 & 0 & 0 \\
        0 & 0 & 0
        \end{pmatrix}
    $$
with respect to $\{\te_1',\te_2',\te_3',\te_4'\}$ and $\{\tf_1',\tf_2',\tf_3'\}$.

\begin{cor}\label{cor: 121206 Cor1}
$\tX_2$ is the union of four $G$-orbits through $(B_1^{(A)},B_2^{(A)},{i^{(A)}},{i^{(A)}}^*)$ for $A=\I,\II,\III,\IV$. And each $G$-orbit is irreducible of dimension $12$, $12$, $11$ and $9$ respectively.
\end{cor}

\proof Let $\tX_2'$ be the union of four $G$-orbits.
Since $\tX_2'\subset \tX_2$, we need to show the opposite inclusion.  By Corollary \ref{cor: 121206 cor}, we have $p(\tX_2)\subset \bigcup_{A=\I,\II} \Sp.{i^{(A)}}= p(\tX_2' )$. Thus $p(\tX_2)= \bigcup_{A=\I,\II} \Sp.{i^{(A)}}= p(\tX_2' )$. On the other hand, $p\inv(p(x))\cong \SL(2)$ for $x\in \tX_2$ by \eqref{eq: Sl2}. Therefore $\tX_2\subset p\inv(p(\tX_2))=p\inv(p(\tX_2'))=\tX_2'$.

By Theorem \ref{th: dim of Sp times O} (2), the explicit value of $\dim G.(B_1^{(A)},B_2^{(A)},{i^{(A)}},{i^{(A)}}^*)=3+ \dim (\Sp\times \rO).{i^{(A)}} $ for each $A$, is computed as above.

The irreducibility comes from the fact that the above $G$-orbits are the $(\SL(2)\times \Sp\times \SO)$-orbits. To see this we observe that $-\Id_W\in \rO$ acts on $\bN$ in the same way as $-\Id_V\in\Sp$.  \qed\vskip.3cm

\begin{lem}\label{lem: 121206 Cor1 and 121213}
Let $x\in G.(B_1^{(A)},B_2^{(A)},{i^{(A)}},{i^{(A)}}^*)$ for $A=\I,\II,\III,\IV$. Then we have the following assertions.

$(1)$ If $A=\I$, $x$ is unstable and $\fsp^x$ is trivial.

$(2)$ If $A=\II$, $x$ is stable (hence $\Sp^x$ is trivial).

$(3)$ If $A=\III,\IV$, $x$ is unstable.

$(4)$ $\tX$ is a reduced complete intersection.
\end{lem}

\proof
(1) It is direct to check $e_1+e_3$ is a common eigenvector of $B_1^{(\I)}$ and $B_2^{(\I)}$. We have ${i^{(\I)}}^*(e_1+e_3)={i^{(\I)}}^*(\te_3)=0$. So $e_1+e_3$ violates the costability of $x$.

We check $\fsp^x$ is trivial. Use \eqref{eq: ft fp}. Then we have
    \begin{equation*}
    \begin{aligned}
    &
    \fsp^{B_1^{(\I)}}
    =\left\{
        \begin{pmatrix} P & JQ^t \\ JQ & S \end{pmatrix}\in \fsp \middle|\, [
        \begin{pmatrix} P & JQ^t \\ JQ & S \end{pmatrix},
        \begin{pmatrix} 0 & I \\ I & 0 \end{pmatrix}] =0
        \right\}
    \\
    &
    \quad \quad \quad
    =\left\{ \begin{pmatrix} P & R \\ R & P \end{pmatrix} \middle|\, \tr P=0, \ R=JQ,\ Q=Q^t \right\},
    \\
    &
    \fsp^{B_1^{(\I)}}
    \cap \fsp^{B_2^{(\I)}}
    = \left\{\begin{pmatrix} P & R \\ R & P \end{pmatrix} \middle|\, [\begin{pmatrix} P & R \\ R & P \end{pmatrix},\begin{pmatrix} I & H \\ -H & -I \end{pmatrix}]=0,\ \tr P=0 \right \}
    \\
    &
    \quad \quad \quad =\cc\left\langle \begin{pmatrix} X & -X \\ -X & X \end{pmatrix}, \begin{pmatrix} Y & Y \\ Y & Y \end{pmatrix}, \begin{pmatrix} H & 0 \\ 0 & H \end{pmatrix}\right\rangle.
    \end{aligned}
    \end{equation*}
On the other hand, $\Image(i^{(\I)})=\cc\langle \te_2,\te_3\rangle=\cc\langle e_1+e_3,e_2-e_4\rangle$. So
    \begin{equation*}
    \begin{aligned}
    &
    \fsp^{i^{(\I)}}
    =\{g\in \fsp|\, g(e_1+e_3)=g(e_2-e_4)=0\}.
    \end{aligned}
    \end{equation*}

We claim $\fsp^{B_1^{(\I)}}\cap \fsp^{B_2^{(\I)}}\cap \fsp^{{i^{(\I)}}}=0$. The solution $a,b,c\in\cc$ of the following equations
    \begin{equation*}
    \begin{aligned}
    &
    \left(a \begin{pmatrix} X & -X \\ -X & X \end{pmatrix} + b \begin{pmatrix} Y & Y \\ Y & Y \end{pmatrix} + c \begin{pmatrix} H & 0 \\ 0 & H \end{pmatrix}\right)(e_1+e_3)=0
    \\
    &
    \left(a \begin{pmatrix} X & -X \\ -X & X \end{pmatrix} + b \begin{pmatrix} Y & Y \\ Y & Y \end{pmatrix} + c \begin{pmatrix} H & 0 \\ 0 & H \end{pmatrix}\right)(e_2-e_4)=0
    \end{aligned}
    \end{equation*}
is trivial. So the claim is proven.

(2) We show there does not exist a common eigenvector $v$ of $B_1^{(\II)}$ and $B_2^{(\II)}$ such that ${i^{(\II)}}^*(v)=0$. By direct computation we have $\Ker(B_1^{(\II)})=\cc\langle e_1,e_3\rangle$ and $\Ker(B_2^{(\II)})=\cc\langle e_1+\sqrt{-1}e_3,e_2+\sqrt{-1}e_4 \rangle$. So we have $\Ker(B_1^{(\II)})\cap\Ker(B_2^{(\II)})=\cc\langle e_1+\sqrt{-1}e_3\rangle$.

On the other hand ${i^{(\II)}}^*(e_1+\sqrt{-1}e_3)=-{i^{(\II)}}^*(\te_2')\neq0$. This proves (2).

(3) In (2) we checked that $e_1+\sqrt{-1}e_3$ is a unique common eigenvector of $B_1^{(A)}$ and $B_2^{(A)}$ for $A=\III,\IV$ up to constant.

On the other hand ${i^{(A)}}^*(e_1+\sqrt{-1}e_3)=-{i^{(A)}}^*(\te_2')=0$ for $A=\III,\IV$. So $e_1+\sqrt{-1}e_3$ violates the costability of $x$.

(4)
By Corollaries \ref{cor: tX1}, \ref{cor: 121206 Cor1} and the dimension reason, $\tX$ is a complete intersection of dimension 12. By (1) and (2), $\tX$ is smooth along the two 12-dimensional $G$-orbits since $\fsp^x=0$ implies that $d\mu_x\colon T_x\bN' \to\fsp$ is surjective. Therefore $\tX$ is reduced along the two 12-dimensional $G$-orbits. Hence $\tX$ itself is reduced by \cite[Prop.\ 5.8.5]{EGA4}. \qed\vskip.3cm


\subsection{Description of $\mu\inv(0)\git \Sp$}
\label{subsec: description of mu inv 0 git Sp}

In the previous subsection we proved
    \begin{enumerate}
    \item $\tX$ is a reduced complete intersection with two irreducible components $\obS$ and $\obU$, where $\bS$ is the stable locus and $\bU$ is a $G$-orbit $G.(B_1^{(\I)},B_2^{(\I)},{i^{(\I)}},{i^{(\I)}}^*)$.
    \item $\bS$ is a $G$-orbit $G.(B_1^{(\II)},B_2^{(\II)},i_4,i_4^*)$.
    \item The $\Sp$-action on $\bU$ is locally free (i.e., $\fsp^x=0$ for any $x\in \bU$).
    \item $\tX\git \Sp=\obS\git \Sp\cup \obU\git\Sp$ as varieties (as $\tX$ is reduced).
    \end{enumerate}

From now on we fix an orthogonal basis of $W$. Then $\fo$ is the set of anti-symmetric matrices.
Let us identify
    \begin{equation}\label{eq: oW cc3}
    \fo \cong \cc^3,\quad \begin{pmatrix} 0 & e & f
    \\ -e & 0 & g \\ -f & -g & 0 \end{pmatrix}\mapsto (e,f,g) .
    \end{equation}
The characteristic polynomial of $A:= \begin{pmatrix} 0 & e & f
    \\ -e & 0 & g \\ -f & -g & 0 \end{pmatrix}$
in $t$ is $t^3+(e^2+f^2+g^2)t$. Therefore $A$ is nilpotent if and only if $e^2+f^2+g^2=0$. Since any nonzero nilpotent element $x$ in $\fo$ has the $a$-diagram $aaa$, the minimal nilpotent orbit is $\rO.x$. Hence, $\bP$ is the quadric surface in $\cc^3$ defined  by  $e^2+f^2+g^2=0$, which also equals the nilpotent variety.

\begin{lem}\label{lem: N}
The map $[(x,y)]\mapsto (x^2,\sqrt{-1}xy,y^2)$ gives an isomorphism $\cc^2/ \zz_2\cong \bP$, and hence $\bP$ is an irreducible normal variety.
Moreover, $\bP^{\rank 2}:=\{ A\in \bP |\, \rank A=2\}=\bP \setminus 0\cong (\cc^2\setminus0) / \zz_2$.
\end{lem}

\proof The first isomorphism is well-known in invariant theory. Since $\cc^2$ is irreducible and normal, so is $\cc^2/ \zz_2$.

We prove the second assertion. The $a$-diagram of a nilpotent matrix $A\in \fo$ is either
    \begin{equation}
    \xy <1cm,0cm>:
    (1,.4)*{a};
    (1,0)*{a}  , (1,-0.4)*{a}
    \endxy
    \quad \mbox{or}\quad
    \xy <1cm,0cm>:
    (1,0)*{aaa}.
    \endxy
    \end{equation}
Therefore $A\neq 0$ means $\rank A=2$.\qed\vskip.3cm

\begin{lem}\label{lem: Phi tX1}
$\Phi_{\tX_1}\colon\tX_1\to \bP$, $(B_1,B_2,i,i^*)\mapsto (\tr(B_1^2),\sqrt{-1}\tr(B_1 B_2),\tr(B_2^2))$, is the GIT quotient by $\Sp$.
\end{lem}

\proof Let $x:=(B_1,B_2,i,i^*)\in \tX_1$. By Corollary \ref{cor: tX1} (2), $x$ is unstable. Suppose $\Sp.x$ is closed in $\bN$. Since $x\notin \mu\inv(0)^\reg$, we have $x^s=0$ by Theorem \ref{th: main 2} and Theorem \ref{th: k2} (4). Thus $i=0$. Let $T:=\{(aB,bB,0,0)|\, B\in \fp',\ a,b\in\cc\}$. Let $\phi:=\Phi_{\tX_1}|_T$. It is enough to show that $\phi$ is the GIT quotient by $\Sp$ by Lemma \ref{lem: 121129 Claim (2)}.

Since the $\Sp$-action on $\bP$ is trivial and $\bP$ is normal, we need to show that $\phi\git\Sp\colon T\git\Sp\to \bP$ is bijective by Zariski's main theorem.

Since $\tr (v_1^2)\neq0$, $\phi$ is surjective and thus so is $\phi\git\Sp$.

To show injectivity, it is enough to show that $\phi\inv(c)$ is an $\Sp$-orbit for any $c\in \bP\setminus0$ and that 0 is the unique closed $\Sp$-orbit in $\phi\inv(0)$.

Let $c\in \bP\setminus0$. Then $c=(a^2,\sqrt{-1}ab,b^2)$ for some $a,b\in \cc$. Then we have $\phi\inv(1,0,0)\cong \phi\inv(c)$, $(B,0,0,0)\mapsto (aB,bB,0,0)$. Thus $\phi\inv(c)$ is an irreducible variety of dimension 4 since $\dim T=6$ and $\dim \bP=2$.

On the other hand, we have ${\fp'}^B=\cc\langle B\rangle$ for any $B\in \fp'\setminus 0$ by Lemma \ref{lem: 121129 Claim (2)}. So $\dim {\fp}^B=2$. By \cite[Prop.\ 5]{KR} we have $\dim \ft^B=6$ and $\dim \Sp.B=4$. This means $\phi\inv(c)$ is a $\Sp$-orbit by irreducibility and the dimension reason.  This proves the first item.

By a similar argument we have $\varphi\git\Sp\colon \fp'\git\Sp\to \cc$ is a birational surjective morphism, where $\varphi\colon \fp'\to \cc$ is given by $B\mapsto \tr B^2$. Since both $\fp'\git\Sp$ and $\cc$ are irreducible normal varieties of dimension 1, $\varphi\git\Sp$ is an isomorphism by Zariski's main theorem. Therefore 0 is the unique closed $\Sp$-orbit in $\varphi\inv(0)$. This proves the second item.
\qed\vskip.3cm

\begin{lem}\label{lem: Phi S}
$\Phi_{\obS}\colon \obS \to \bP$, $(B_1,B_2,i,i^*)\mapsto i^*i$, is the GIT quotient by $\Sp$.
\end{lem}

\proof Let $x:=(B_1,B_2,i,i^*)\in \bS$. By the $ab$-diagram \eqref{eq: ab III} of $i$, we have $i^*i\in \bP^{\rank2}$. So $\Phi_{\obS}(\bS)=\bP^{\rank2}$ and $\Phi_{\obS}$ is well-defined.

We claim that $\Phi_{\obS}\inv(a)$ is an $\Sp$-orbit for any $a\in \bP\setminus0=\bP^{\rank2}$. Note that since  $\bS$ is irreducible, so is $\Phi_{\obS}\inv(a)$. Since $x$ is stable, $\Sp^x$ is trivial. Thus $\dim \Sp.x=10$. On the other hand, $\Phi_{\obS}\inv(i^*i) = (\SL(2)\times \Sp).x$ by Theorem \ref{th: first fund th Inv th} and \eqref{eq: Sl2}. Thus $\Phi_{\obS}\inv(i^*i)$ is an irreducible 10-dimensional variety.  If $\Sp.x\subsetneq(\SL(2)\times \Sp).x$ then for $y\in  (\SL(2)\times \Sp).x\setminus \Sp.x$, we have $\dim \Sp.y=10$. But then $\Phi_{\obS}\inv(i^*i)$ contains two disjoint locally closed subvarieties $\Sp.x$ and $\Sp.y$ of dimension 10. This contradicts the irreducibility of $\Phi_{\obS}\inv(i^*i)$. Therefore $\Sp.x=\Phi_{\obS}\inv(i^*i)$ as desired.

Since $\bS$ consists of $\Sp$-closed orbits (Theorem \ref{th: main 2}), $\bS/\Sp$ is Zariski open in $\obS\git \Sp$. By Luna's slice theorem \cite{Lu}, $\bS/\Sp$ is a smooth variety. By the above claim and Zariski's main theorem, $\oPhi_{\obS}|_{\bS/\Sp}\colon \bS/\Sp\to \bP^{\rank2}$ is an isomorphism, where $\oPhi_{\obS}\colon \obS\git\Sp\to \bP$ is the induced morphism. Let us finish the proof of the lemma. Let $f\in \cc[\obS]^{\Sp}$. Then $f|_S\in \Phi^*\Gamma(\cO_{\bP\setminus0})$. By the normality of $\bP$, $\Gamma(\cO_{\bP\setminus0})= \cc[\bP]$. Thus $f\in \Phi_{\obS}^*\cc[\bP]$. This means $\cc[\obS]^{\Sp}=\Phi_{\obS}^*\cc[\bP]$, equivalently $\Phi$ is the GIT quotient by $\Sp$ as desired.
\qed\vskip.3cm

\begin{lem}\label{lem: tX1 obU}
$\tX_1\subset \obU$ and $\tX_1\git \Sp=\obU\git\Sp$.
\end{lem}

\proof Since $\bU$ is an irreducible reduced locally free $G$-orbit of dimension 12, we have $\obU\git \Sp$ is an irreducible reduced variety of dimension $\le 2$.
Now the second assertion follows from the first by Lemma \ref{lem: Phi tX1}.

Let us prove the first assertion. Suppose $\tX_1\nsubseteq\obU$. Then $\tX_1\subset \obS$ since $\tX_1$ is irreducible by Corollary \ref{cor: tX1} (1). Let $T:=\{(aB,bB)\in \fp'\times\fp'|\, a,b\in\cc,\  B\in \fp'\}$. By Lemma \ref{lem: 121129 Claim (1)} and \eqref{eq: ab IV}, $\tX_1$ is the closure of $T\times (\Sp\times \rO).i$ for a nonzero $i\in \Hom(W,V)$ with $i^*i=0$. By Lemma \ref{lem: Phi S}, $\Phi_{\obS}(\tX_1)=0$. This contradicts Lemma \ref{lem: Phi tX1}. \qed\vskip.3cm

\begin{defn}
Define $\Phi\colon \tX\to (\bP\times 0)\cup (0\times \bP)\ (\subset \bP\times \bP)$, $(B_1,B_2,i,i^*)\mapsto ((\tr(B_1^2),\sqrt{-1}\tr(B_1 B_2),\tr(B_2^2)),i^*i)$.
\end{defn}

To see $\Phi$ is well-defined morphism one notices that ${i^{(\I)}}^*{i^{(\I)}}={i^{(\II)}}^*{i^{(\II)}}={i^{(\III)}}^*{i^{(\III)}}=0$ and $\tr({B_1^{(\II)}}^2)=\sqrt{-1}\tr(B_1^{(\II)} B_2^{(\II)})=\tr({B_2^{(\II)}}^2)=0$, which come from Corollary \ref{cor: 121206 cor} and the direct computation respectively.

\begin{thm}\label{th: 121227 Claim 1}
$\Phi$ is the GIT quotient by $\Sp$ onto $(\bP\times 0)\cup (0\times \bP)$. Hence, $\mu\inv(0)\git \Sp\cong \cc^2\times \left((\bP\times 0)\cup (0\times \bP)\right)$. \end{thm}

\proof  By Lemmas \ref{lem: Phi tX1}, \ref{lem: Phi S} and \ref{lem: tX1 obU}, $\Phi$ factors through $\tX\git \Sp$. The canonical morphism given by the composite
    $$
    \xymatrix{
    (\bP\times 0)\cup (0\times \bP) \ar[r]^{f}  & \tX_1\git\Sp \coprod_{[0]} \obS\git \Sp\ar[r] & \tX\git \Sp
    }
    $$
is the inverse of $\Phi\git \Sp$, where $f:=(\Phi|_{\tX_1}\git \Sp)\inv \coprod_{[0]} (\Phi|_{\obS}\git \Sp)\inv $ and $\Phi\git\Sp\colon \tX\git \Sp\to (\bP\times 0)\times(0\times \bP)$ is the induced morphism. \qed\vskip.3cm



\section{Moduli spaces of $\SO(2)$-data with $k=4$}
\label{sec: proof of Theorem N2 (2)}

This section will be devoted to the proof of Theorem \ref{th: N2} (2).

Let $\dim V=k=4$ and $\dim W=N=2$.
Let $\fp':=\{B\in \fp|\, \tr(B)=0\}$. Let $\tX:=\{(B_1,B_2,i,i^*)\in \mu\inv(0)|\, B_1,B_2\in \fp'\}$. Then as in \S\ref{subsec: description of mu inv 0 N3k4},
$\mu\inv(0)\cong \cc^2\times \tX$. Let $p\colon \tX\to \Hom(W,V)$ be the projection.

\begin{lem}\label{lem: ab for V4W2}
Let $i\in \Image(p)$. Then $ii^*$ is nilpotent and the $ab$-diagram of $i$ is one of the followings:
    \begin{equation}
    \label{eq: 99}
    {\xy <1cm,0cm>:
    (1,0.8)*{ab };
    (1,0.4)*{ba };
    (1,0 )*{b  };
    (1,-0.4)*{ b };
    \endxy} \quad\quad\quad
    {\xy <1cm,0cm>:
    (1,1.2)*{b };
    (1,0.8)*{b  };
    (1, 0.4)*{  b };
    (1, 0)*{  b };
    (1, -0.4)*{  a };
    (1, -0.8)*{  a };
    \endxy}\quad\quad\quad
    {\xy <1cm,0cm>:
    (1,0.4)*{bab };
    (1,0)*{bab };
    \endxy} \quad\quad\quad
    {\xy <1cm,0cm>:
    (1,0.8)*{bab };
    (1,0.4)*{b };
    (1, 0 )*{ b };
    (1, -0.4 )*{ a };
    \endxy}
    \end{equation}
\end{lem}

\proof Let $(B_1,B_2,i,i^*)\in \tX$. By Lemma \ref{lem: 121206 claim pre}, $([B_1,B_2])^2=(ii^*)^2$ is a scalar. Since $\rank \, i\le2$ and $\dim V=4$, the scalar is 0. Since there is no nontrivial nilpotent element in $\fo$, we have $i^*i=0$. By a similar argument as in \S\ref{subsec: proof Theorem k2}, the $ab$-diagram of $i$ is one of \eqref{eq: 99} from Table \ref{table: alpha beta}. \qed\vskip.3cm

Let $\tX_1:=\{(B_1,B_2)\in \fp'\times\fp'|\, [B_1,B_2]=0\}\times \{i\in \Hom(W,V)|\, ii^*=0\}$. Let $\tX_2:=\tX\setminus \tX_1$.

\begin{thm}
$\mu\inv(0)^\reg=\emptyset$.
\end{thm}

\proof
Let $x:=(B_1,B_2,i,i^*)\in \tX$. We will prove that $x$ is not costable.

Suppose $x\in\tX_1$. The $ab$-diagram of $i$ is either the first or the second in \eqref{eq: 99}. Therefore $\dim\Ker(i^*)\ge3$. As in the proof of Corollary \ref{cor: tX1} (2), we see that $x$ is not costable.

Suppose $x\in \tX_2$. As in the proofs of Lemma \ref{lem: 121206 Cor1 and 121213} (1) and (3), we deduce that $x$ is not costable.  \qed
\vskip.3cm


\appendix


\section{Finite dimensionality of weight spaces}\label{sec: finite dim}

The main purpose of this section is to prove that each weight space of $\cc[\mu\inv(0)]^{G(V)}$ with respect to $T$ is finite-dimensional.

Let $\ep(m)=0$ (resp.\ $\ep(m)=1$) if $m$ is even integer (resp.\ odd integer). If $\vep=-1$ then using a symplectic basis of $V$, we identify $V=\cc^k$. If $\vep=+1$ then using the basis
    $$
    \begin{aligned}
    \{
    &
    f_1\pm\sqrt{-1}f_2,f_3\pm\sqrt{-1}f_4,...,f_{k-\ep(k)-1}\pm\sqrt{-1}f_{k-\ep(k)}, (f_k)\}
    \end{aligned}
    $$
we identify $V=\cc^k$, where $\{f_1,f_2,...,f_k\}$ is an orthogonal basis of $V$. Here the notation $(f_k)$ denotes $f_k$ only when $k$ is odd (vacuous otherwise). By a similar way we identify $W=\cc^N$. Let $\lfloor a \rfloor$ be the maximal integer in $\zz_{\le a}$, where $a\in \rr$.
We fix maximal tori of $G(V)$ and $G(W)$ as
    $$
    \begin{aligned}
    &
    T_{G(V)}=\{\diag(z_1,z_2,...,z_{\lfloor k/2 \rfloor})\oplus \diag(z_1\inv,z_2\inv,...,z_{\lfloor k/2 \rfloor}\inv)|\, z_1,z_2,...,z_{\lfloor k/2 \rfloor}\in \cc^*\}, \\
    &
    T_{G(W)}=\{\diag(t_1,t_2,...,t_{\lfloor N/2 \rfloor})\oplus \diag(t_1\inv,t_2\inv,...,t_{\lfloor N/2 \rfloor}\inv)|\, t_1,t_2,...,t_{\lfloor N/2 \rfloor}\in \cc^*\}
    \end{aligned}
    $$
respectively.

Now we identify the rings of characters of $T_{G(V)}$, $T_{G(W)}$ and $(\cc^*)^2$
    $$
    \begin{aligned}
    &
    R(T_{G(V)})=\zz[z_1^{\pm1},z_2^{\pm1},...,z_{\lfloor k/2 \rfloor}^{\pm1}],
    \\
    &
    R(T_{G(W)})=\zz[t_1^{\pm1},t_2^{\pm1},...,t_{\lfloor N/2 \rfloor}^{\pm1}]
    \\
    &
    R((\cc^*)^2)=\zz[q_1^{\pm 1},q_2^{\pm1}]
    \end{aligned}
    $$
respectively.

Let us prove
    \begin{equation}\label{eq: cX mu}
    \cX_{\cc[\mu\inv(0)]^{G(V)}}\in \hR(T):=R(T_{G(W)})[[q_1\inv,q_2\inv]]
    \end{equation}
where $T=T_{G(W)}\times (\cc^*)^2$ and $\cX_{\cc[\mu\inv(0)]^{G(V)}}$ is the formal $T$-character.
The idea is to use the following two:
    \begin{enumerate}
    \item
    $\mu\inv(0)\git G(V)$ is a closed $(\cc^*)^2\times G(W)$-subscheme of $\bN\git G(V)$;
    \item
    the (surjective) GIT quotient $\bN\to \bN\git G(V)$ is $(\cc^*)^2\times G(W)$-equivariant.
    \end{enumerate}
So each weight space of $\cc[\mu\inv(0)]^{G(V)}$ is a subquotient of the weight space of the same weight of $\cc[\bN]$. So the proof will be done if we show a stronger claim:
    $$
    \cX_{\cc[\bN]}\in \hR(T).
    $$
The claim will follow from an even stronger claim:
    $$
    \cX^{(\cc^*)^2}_{\cc[\bN]}\in \zz[[q_1^{-\frac12},q_2^{-\frac12}]]
    $$
where $\cX^{(\cc^*)^2}_{\cc[\bN]}$ denotes the formal $(\cc^*)^2$-character of $\cc[\bN]$. This is because all the monomials in $\cX_{\cc[\bN]}$ have nonnegative integer coefficients.

Let us check the last claim. The decomposition $\bN=\fp(V)\oplus \fp(V)\oplus \Hom(W,V)$ is the weight decomposition with respect to $(\cc^*)^2$.
The first two direct summands $\fp(V)$ of $\bN$ are of weight $q_1$ and $q_2$ respectively. The last summand $\Hom(W,V)$ is of weight $(q_1q_2)^{\frac12}$. Since $\cc[\bN]=S(\fp(V)\dual)\otimes S(\fp(V)\dual)\otimes S(\Hom(W,V)\dual)$, where $S$ denotes the symmetric product,
we have
    $$
    \cX^{(\cc^*)^2}_{\cc[\bN]}=\left(\sum_{n\ge0}q_1^{-n}\right)^{\frac12 k(k+\vep)} \left(\sum_{n\ge0}q_2^{-n}\right)^{\frac12 k(k+\vep)} \left(\sum_{n\ge0}(q_1q_2)^{-\frac n2}\right)^{kN}.
    $$
This proves the last claim.

To complete the proof of \eqref{eq: cX mu} we recall from \S\ref{intro: subsec: ADHM} that $i\in \Hom(W,V)$ always appears together with $i^*$ in $\cc[\bN]^{G(W)}$. Hence any monomials with non-integer exponents in $\cX_{\cc[\mu\inv(0)]^{G(W)}}$ have coefficient 0. This proves \eqref{eq: cX mu}.

\begin{rk}
(1) $\hR(T)$ is a ring.

(2) The origin 0 is the unique $(\cc^*)^2$-fixed point of $\bN$.
\end{rk}


\section{Scheme structure of moduli space of framed vector bundles with symplectic and orthogonal structure}\label{sec: moduli space of G bundles}

Let $V,W$ be vector spaces of dimension $k,N$ with $(\,,\,)_\vep,(\,,\,)_{-\vep}$ respectively, where $\vep=\pm1$.
The above pairings give isomorphisms $a_\vep\colon V\to V\dual$ and $b_{-\vep}\colon W\to W\dual$ given by $v\mapsto (v,\bullet)_\vep$ and $w\mapsto (w,\bullet)_{-\vep}$ respectively.
Then, $a_\vep^\vee =\vep a_\vep $ and $b_{-\vep}\dual=-\vep b_{-\vep}$.
We use the notations $\bM$ and $*_\bM$ in \S\ref{subsec: the Orbit-closedness and the semisimplicity}.
Let $\ast:=*_\bM$ and $x^*:=\ast(x)$ for simplicity.
Let $\mu_\bM\colon \bM\to \gl(V)$ be the moment map given by $(B_1,B_2,i,j)\mapsto [B_1,B_2]+ij$.
It is obvious that $\mu=\mu_\bM|_\bN$.

We define an involution $\oast\colon \bM^\reg/ \GL(V)\to \bM^\reg/ \GL(V)$ by $\GL(V).x\mapsto \GL(V).x^*$.
Then the fixed locus $\bM^\reg/\GL(V)^\oast$ (resp.\ $\mu_\bM\inv(0)^\reg/\GL(V)^\oast$) is a smooth subscheme of $\bM^\reg/\GL(V)$ (resp.\ $\mu_\bM\inv(0)^\reg/\GL(V)$).

Since $\bN$ is the fixed locus of $\bM$, we have a canonical embedding $\iota\colon\mu\inv(0)^\reg/G(V)\to \mu_\bM\inv(0)^\reg/\GL(V)^\oast$.
We will see $\iota$ is surjective.

Let $\cG$ be the vector bundle locus of the Gieseker moduli scheme of framed torsion-free sheaves $E$ with rank $N$ and $c_2(E)=k$.
By Barth's correspondence \cite{Barth}, there is an isomorphism $F'\colon \mu_\bM\inv(0)^\reg/\GL(V)\to \cG$ (cf.\ \cite[\S2]{Lecture}).
We denote by the same notation $\oast$ the induced involution on $\cG$ via $F'$.
Thus we have the isomorphism between the fixed loci $F\colon \mu_\bM\inv(0)^\reg/\GL(V)^\oast\to \cG^\oast$ by restriction.

Donaldson's argument \cite{Do} asserts that the image of the composite $F\iota$ is the set of the isomorphism classes of framed vector bundles $E$ which admits an isomorphism $\phi\colon E\to E\dual$ with $\phi\dual=-\vep\phi$ and $\phi|_x=b_{-\vep}$ (after identifying $E|_x=W$ via the given framing for any $x\in l_\infty$).
See also \cite{BrSan}.

We claim that the image of $F$ is contained in that of $F\iota$ given as above.
So our claim will assert that $\iota$ is surjective.
The proof itself goes along Donaldson's argument, so will be sketchy.
Let $x\in \mu_\bM\inv(0)^\reg$ such that $g.x=x^*$ for some $g\in \GL(V)$.
By taking $*$ to the both sides of $g.x=x^*$, we obtain $g.x=g^*.x$.
By the stability of $x$, we have $g=g^*$.
The cohomology sheaf of the monad associated to $x^*$ is isomorphic to $E\dual$.
The induced maps by $g,g^*$ between the monads are explicitly written in terms of linear maps in $\End(V^{\oplus2}\oplus W)$.
By diagram-chasing and passing to the cohomology sheaves, the constraint $g=g^*$ gives an isomorphism $\og\colon E\to E\dual$ with $\og\dual=-\vep\og$ and $\og|_x=b_{-\vep}$.
This finishes the proof of the claim.

By Zariski's main theorem we obtain the isomorphism $\cM^K_n\cong\cG^\oast$, which was used in \S\ref{intro: subsec: ADHM2}.


\section{Proof of Proposition \ref{prop: smooth} (2)}\label{app: B}

We prove Proposition \ref{prop: smooth} (2) in this section: the codimension of
$m\inv(X)\setminus (\pi_1\inv(\fp_k)\cup \pi_2\inv(\fp_k))$ in $m\inv (X)$ is larger than 1.

\begin{lem}\label{lem: one-orbit}
$($\cite[Theorem\ XI.4]{Ga}$)$
Let $B\in \fp$. Then $\rO(V).B=\fp\cap \GL(V).B $, where $\rO(V).B$ and $\GL(V).B$ are the orbits by conjugation. \qed \end{lem}

Let $E_{\gl}\colon\gl\to S^k\cc$ be the morphism mapping $B$ to the unordered set of eigenvalues of $B$. Here, $S_k$ is the symmetric group of $k$-letters acting on $\cc^k$ by permutation of coordinates, so that $S^k\cc:=\cc^k/S_k$. Let $E:=E_{\gl}|_{\fp}$ and $E_l:=E|_{\fp_l}$. To construct $E_\gl$ explicitly, let $P\colon S^k\cc\to \cc^k$ be the isomorphism by $[(a_1,...,a_k)]\mapsto (p_1(a),...,p_k(a))$, where $a:=(a_1,...,a_k)$ and $p_i(a)=a_1^i+...+a_k^i$ (the $i^{\mathrm{th}}$ power sum).
Let $E'_{\gl}\colon\gl\to \cc^k$ by $ B\mapsto (\tr B,\tr B^2,...,\tr B^k) $. Let $E_{\gl}:=P\inv\circ E'_\gl$.

Let $\fp_l^{(e)}:=\{B\in \fp_l|\, \mbox{$B$ has $e$ distinct eigenvalues}\}.$ Let $\fp_l^{(\le e)}:=\bigsqcup_{e'\le e} \fp_l^{(e')}$. Then
$\fp_l^{(\le e)}$ is a closed subvariety of $\fp$. Indeed, let $\Delta^{(e)}\subset S^k\cc$ be the locus of all the unordered sets of $e$ distinct points. Let $\Delta^{(\le e)}:=\bigsqcup_{e'\le e} \Delta^{(e')}$. Then  $\Delta^{(\le e)}$ is a closed subvariety of $S^k\cc$. Therefore $E\inv(\Delta^{(\le e)})$ is a closed subvariety of $\fp$. In particular, $\fp_l^{(e)}=\fp_l\cap E\inv(\Delta^{(e)})$ is locally closed in $\fp$ by Lemma \ref{lem: local closed}. It is manifest that if $e>k$ then for any $l$, $\fp_l^{(e)}=\emptyset$ and that if $e=k$ then $\fp_l^{(e)}\neq\emptyset$ if and only if $l=k$.

\begin{lem}\label{lem: E}
$(1)$ If $\fp_l^{(e)}\neq\emptyset$ then  $\dim\fp_l^{(e)}\le \dim \fp-l+e$. In particular, if $l>k$ then $\dim\fp_l^{(e)}\le \dim \fp-2$.

$(2)$ If $e=k-1$ then $\fp_l^{(e)}\neq \emptyset$ if and only if $l\in\{k,k+1\}$.

$(3)$ $\fp_{k+1}^{(k-1)}$ consists of $B\in \fp$ conjugate by $\rO(V)$ to
        $  \diag (a_1,a_1,a_2,...,a_{k-1})$, where $a_1,a_2,...,a_{k-1}$ are distinct in $\cc$.
\end{lem}

\proof (1) The image $E_l|\colon\fp_l^{(e)}\to S^k\cc$ is contained in $\Delta^{(e)}$. Any nonempty fibre of $E_l$ is a union of $\rO(V).B$ for finitely many $B\in \fp_l^{(e)}$ (Lemma \ref{lem: one-orbit}), so that its dimension is $\dim \rO(V).B=\dim \rO(V)-\dim \rO(V)^B= \dim \ft-\dim \ft^B= \dim \fp-\fp^B=\dim\fp-l$, where the third identity comes from Lemma \ref{lem: KR}. Therefore $\dim \fp_l^{(e)}\le \dim\fp-l+\dim \Delta^{(e)}$. Since $\dim\Delta^{(e)}=e$, we have proven (1).

(2)
Let $e=k-1$. Let $B\in \fp_l^{(e)}$. Let $a_1,...,a_{k-1}$ be the (distinct) eigenvalues of $B$. We may assume that only the $a_1$-eigenspace of $B$ is 2-dimensional while the other ones are all 1-dimensional. The Jordan normal form of $B$ is either $\begin{pmatrix} a_1 & 1 \\ 0 & a_1\end{pmatrix}\oplus \diag(a_2,...,a_{k-1})$ or $\diag(a_1,a_1,...,a_{k-1})$.  Both cases actually happen, since $\begin{pmatrix} a_1 & 1 \\ 0 & a_1\end{pmatrix}$ is conjugate by $\GL(2)$ to a symmetric matrix $\begin{pmatrix} a_1 +\sqrt{-1}  & 1 \\ 1 & a_1 -\sqrt{-1}\end{pmatrix}$.  We have $\dim \gl^B= k$ and $k+2$ respectively. Using $\gl^B=\fp^B\oplus \ft^B$ (by \eqref{eq: ftp}) and Lemma \ref{lem: KR}, we have $\dim \fp^B= k$ and $k+1$ respectively. This proves (2).

(3) follows from Lemma \ref{lem: one-orbit}. \qed\vskip.3cm

\begin{lem}\label{lem: F}
Let $X\in \ft$.  Let $i\in \{1,2\}$.
Suppose $\pi_i\inv(\fp_l^{(e)})\cap m\inv(X)\neq\emptyset$. Then $\dim \pi_i\inv(\fp_l^{(e)})\cap m\inv(X)\le \dim \fp +e$. In particular, if $e\le k-2$ then $\dim \pi_i\inv(\fp_l^{(e)})\cap m\inv(X)\le \dim \fp+k -2=\dim m\inv(X)-2$. \end{lem}

\proof  We claim that any nonempty fibre of $\pi_i|\colon\pi_i\inv(\fp_l^{(e)})\cap m\inv(X)\to \fp_l^{(e)}$ is of dimension $l$. Take $B\in \fp^{(e)}$ and identify $\pi_i\inv(B)$ with $\fp$. Then $\pi_i\inv(B)\cap m\inv(X)$, unless empty, is an affine space isomorphic to $\fp^B$, since for any $B'\in \pi_i\inv(B)\cap m\inv(X)$, $B'-B\in \fp^B$. The base dimension $\dim \fp_l^{(e)}$ is estimated in Lemma \ref{lem: E}. Thus the lemma is proven.  \qed\vskip.3cm

Now we are ready to estimate the codimension of
$m\inv(X)\setminus (\pi_1\inv(\fp_k)\cup \pi_2\inv(\fp_k))$ in $m\inv (X)$. By Lemma \ref{lem: F} and Lemma \ref{lem: E} (2), to check the codimension $\ge 2$, it suffices to check that so is the codimension of $ \pi_i\inv(\fp_{k+1}^{(k-1)})\cap \pi_j\inv(\fp\setminus \fp_k)\cap m\inv(X)$ in  $m\inv(X)$, whenever $\{i,j\}=\{1,2\}$. By Lemma \ref{lem: F}, $\pi_i\inv(\fp_{k+1}^{(k-1)}) \cap m\inv(X)$ is of codimension $\ge1$ in $m\inv(X)$. It remains to prove $ \pi_i\inv(\fp_{k+1}^{(k-1)})\cap \pi_j\inv(\fp_k)\cap m\inv(X)$ is Zariski (open) dense in $\pi_i\inv(\fp_{k+1}^{(k-1)}) \cap m\inv(X)$. Let $B_1\in \fp^{(k-1)}_{k+1}$. This is reduced to check that
    \begin{equation}\label{eq: star}
    \pi_i\inv(B_1)\cap \pi_j\inv(\fp_k)\cap m\inv(X)\neq \emptyset \ \mbox{provided}\ \pi_i\inv(B_1)\cap m\inv(X)\neq\emptyset
    \end{equation}
since $\pi_i\inv(B_1)\cap m\inv(X)\cong \fp^{B_1}$ irreducible (see the proof of Lemma \ref{lem: F}). Let $B_0\in \pi_i\inv(B_1)\cap m\inv(X)\subset \fp$, where $\pi_i\inv(B_1)$ are canonically identified with $\fp$. Let us write $B_1=g.\diag(a_1,a_1,a_2,...,a_{k-1})$, where $g\in \rO(V)$ and $a_1,...,a_{k-1}$ are distinct (Lemma \ref{lem: one-orbit} and Lemma \ref{lem: E} (3)). Let $B_2:=g.\diag(b_1,...,b_k)$, where $b_1,...,b_k$ are distinct so that $B_2\in \fp_k$.
By the Zariski openness of $\fp_k$ in $\fp$, there exists $u\in \cc\setminus\{1\}$ such that $(1-u)B_0+uB_2\in \fp_k$ since for $u=1$, $B_2\in \fp_k$. Therefore $B_0+\frac{u}{1-u}B_2\in \fp_k$. Now we have $(B_1,B_0+\frac{u}{1-u}B_2)$ or $(B_0+\frac{u}{1-u}B_2,B_1)\in  \pi_i\inv(B_1)\cap \pi_j\inv(\fp_k)\cap m\inv(X)$, which shows \eqref{eq: star}.
This completes the proof of Proposition \ref{prop: smooth} (2).

\end{document}




\title[Moduli spaces of framed symplectic and orthogonal bundles: Corrigendum and addendum]{Corrigendum and addendum to ``Moduli spaces of framed symplectic and orthogonal bundles on $\pp^2$ and the K-theoretic Nekrasov partition functions''}

\author{Jaeyoo Choy}

\address{Dept. Math., Kyungpook Nat'l Univ., Sangyuk-dong, Buk-gu, Daegu 702-701, Korea}

\email{choy@knu.ac.kr}


 \keywords{moduli spaces, framed symplectic and orthogonal bundles, instantons, K-theoretic Nekrasov partition functions}

\begin{abstract}
We give a scheme-theoretic description of $\mu\inv(0)$ and $\mu\inv(0)\git\Sp(2)$, where $\mu$ is the moment map defining the space of ADHM data associated to framed $\SO(3,\rr)$-instantons with instanton number $1$. 
The irreducible component in $\mu\inv(0)\git\Sp(2)$, other than the Zariski closure of the instanton moduli space, is non-reduced, which is a corrigendum to \cite[Theorem 4.3]{Choy1}.
We further describe this component using the factorization property \cite{BFG}.
\end{abstract}

\maketitle
\setcounter{tocdepth}{1}

\thispagestyle{empty} \markboth{Jaeyoo Choy}{Moduli spaces of framed symplectic and orthogonal bundles: Corrigendum and addendum}


\section{Corrigendum}
\label{sec: Erratum and correction}

The Hamiltonian reduction $\mu\inv(0)\git\Sp(2)$ associated to the framed $\SO(3,\rr)$-instanton moduli space with instanton number $4$ is \textit{not} reduced.
Hence the statement of \cite[Theorem 4.3]{Choy1} should be corrected as:
\begin{theorem*}
Let $N=3$ and $k=4$.
Then we have the following properties.

$(1)$ $\mu\inv(0)$ is a non-reduced complete intersection and a union of two irreducible components.
Hence $\mu$ is flat.

$(2)$ One irreducible component of $\mu\inv(0)$ is the Zariski closure of $\mu\inv(0)^{\reg}$.

$(3)$ The reduced scheme $(\mu\inv(0)\git\Sp(V))_{\mathrm{red}}$ is isomorphic to $\cc^{2}\times(\bP\sqcup_{0}\bP)$.
Moreover the isomorphism restricts to $\mu\inv(0)^{\reg}/\Sp(V)\cong\cc^{2}\times(\bP\setminus0)$.
\end{theorem*}

The (corrected) proof will appear in Section \ref{subsec: corrected proof} after correction of \cite[Lemma 8.10 (1), (4)]{Choy1}.
In Sections \ref{sec: addendum}, \ref{sec: scheme structure of GIT quotient}, we describe further the non-reduced locus of $\mu\inv(0)$ and its GIT quotient.


\subsection{Correction of \cite[Lemma 8.10 (1), (4)]{Choy1}}\label{subsec: corrected lemma}

This erratum is due to the author's wrong computation of the dimension of the stabilizers in \cite[Lemma 8.10 (1)]{Choy1}.
This statement is corrected as follows:  
If $A=\I$, then $x$ is unstable and $\dim\fsp^{x}=2$.
We give the correct computation of $\fsp^{x}$, where $x=(B_{1}^{(\I)},B_{1}^{(\I)},i^{(\I)},{i^{(\I)}}^{*})$.
Using the matrix expression of $\ft=\fsp$ (\cite[(8.3)]{Choy1}), we have
     \begin{equation*}
    \begin{aligned}
    &
    \fsp^{B_1^{(\I)}}
    =\left\{
        \begin{pmatrix} P & JQ^t \\ JQ & S \end{pmatrix}\in \fsp \middle|\, [
        \begin{pmatrix} P & JQ^t \\ JQ & S \end{pmatrix},
        \begin{pmatrix} 0 & I \\ I & 0 \end{pmatrix}] =0
        \right\}
    \\
    &
    \quad \quad \quad
    =\left\{ \begin{pmatrix} P & R \\ R & P \end{pmatrix} \middle|\, \tr P=0, \ R=JQ,\ Q=Q^t \right\},
    \\
    &
    \fsp^{B_1^{(\I)}}
    \cap \fsp^{B_2^{(\I)}}
    = \left\{\begin{pmatrix} P & R \\ R & P \end{pmatrix} \middle|\, [\begin{pmatrix} P & R \\ R & P \end{pmatrix},\begin{pmatrix} I & -H \\ H & -I \end{pmatrix}]=0,\ \tr P=0 \right \}
    \\
    &
    \quad \quad \quad =\cc\left\langle \begin{pmatrix} X & X \\ X & X \end{pmatrix}, \begin{pmatrix} Y & -Y \\ -Y & Y \end{pmatrix}\right\rangle,
    \end{aligned}
    \end{equation*}
where $\{H,X,Y\}$ is the standard $\fsl_{2}$-triple.
(The original computation in the proof of \cite[Lemma 8.10 (1)]{Choy1} is wrong because the  matrices $H,-H$ in the above were interchanged mistakenly).
On the other hand $\Image(i^{(\I)})=\cc\langle \te_2,\te_3\rangle=\cc\langle e_1+e_3,e_2-e_4\rangle$ and this is contained in the kernel of any endomorphisms in the above $\fsp^{B_1^{(\I)}}\cap \fsp^{B_2^{(\I)}}$.
Thus we have
	$$
	\fsp^{x}=\fsp^{B_1^{(\I)}}\cap \fsp^{B_2^{(\I)}}=\cc\left\langle \begin{pmatrix} X & X \\ X & X \end{pmatrix}, \begin{pmatrix} Y & -Y \\ -Y & Y \end{pmatrix}\right\rangle.
	$$
(See Remark \ref{subsec: remark} for $\Sp^{x}$).

The statement of \cite[Lemma 8.10 (4)]{Choy1} should be corrected as follows:
$\tX =\{(B_{1},B_{2},i,i^{*})\mid \tr(B_{1})=\tr(B_{2})=0\} $ is a non-reduced complete intersection with two irreducible components. 
Moreover the Zariski closure of $G.x$ is the only non-reduced component, where $G=\SL(2)\times\Sp\times\rO$ and $x$ is of type I as in the above computation.

The verification is as follows:
Since $\fsp^x\neq0$, the differential $d\mu_x$ of the moment map $\mu$ at $x$ is not surjective, which means that $\tX$ has non-reduced scheme structure along the stratum $G.x$ (cf.\ the proof of \cite[Theorem 3.3]{Choy1}).


\subsection{Proof of Theorem}\label{subsec: corrected proof}

Since $\fsp^{x}$ is not trivial, we need to correct every relevant statement in \cite[Section 8.2]{Choy1}:  $\tX$ should be replaced by $(\tX)_{\mathrm{red}}$, and the $G$-action on $\bU=G.x$ is not locally free (\cite[p.300, the item (3)]{Choy1} and \cite[proof of Lemma 8.14]{Choy1}).

Theorem immediately follows from \cite[Lemma 8.14]{Choy1} once we assure $\dim\obU\git\Sp\le 2$ in the proof of the lemma (without assuming the local freeness of the $G$-action on $\bU$).
We use the decomposition $\mu\inv(0)\git \Sp=\mu\inv(0)^{\reg}/\Sp\sqcup S^{2}\aaa^{2}$ (\cite[Theorem 2.6 (1)]{Choy1}).
Here we used the fact that there is no $\SO(3,\rr)$-instanton with instanton number $2$, i.e.\ $\mu\inv_{k'}(0)=\emptyset$, $k'=2$.
Thus the image of $\bU$ in the GIT quotient is contained in the stratum $S^{2}\aaa^{2}$, which implies $\dim\obU\git\Sp\le 2$.
This completes the proof of Theorem.


\section{Further description on the scheme structure of the non-reduced locus of $\tX$}
\label{sec: addendum}

We describe further details of the scheme structure of $\tX$ along the orbit $G.x$, where $x$ is given in the previous section.
Let
    $$
    \tX_0:=\left\{(B_1,B_2,i)\in\tX\middle| \mbox{$B_1$ has two distinct eigenvalues}\right\} 
    $$
(an affine open subset of $\tX$).
Here we identified $\tX$ with the image of a projection so that $\tX$ consists of triples rather than quadruples.
Let $\rho\colon\Hom(W,\cc^2)\to \fsp(1),\ i\mapsto ii^*$, where $W\,(=\cc^3),\cc^2$ are the orthogonal and symplectic vector spaces respectively.
We used this map in the proof of \cite[Theorem 4.1]{Choy1}.

\begin{prop}\label{prop: smooth morphism}
$(1)$
There is a surjective smooth morphism from $\Sp\times\cc^*\times\cc\times (\rho\inv(0))^2$ to $\tX_0$.

$(2)$
This morphism restricts to a smooth morphism from the open subset $\Sp\times\cc^*\times\cc\times (\rho\inv(0)\setminus0)^2$ to $G.x$.
\end{prop}

The proof  and also the precise definition of the above smooth morphism will appear in Section \ref{subsec: proof of prop}.

As a corollary, the structures of local analytic singularities of $\tX_0$ and $(\rho\inv(0))^2$ (resp.\ $G.x$ and $(\rho\inv(0)\setminus0)^2$) are equivalent.
By \cite[Theorem 4.1 (4)]{Choy1}, $\rho\inv(0)$ and $\rho\inv(0)\setminus0$ are non-reduced and thus so are $\tX_0$ and $G.x$.


\subsection{Strata of $\tX_{0}$}

Recall the decomposition $\tX=\tX_{1}\sqcup \tX_{2}$, where $\tX_{1}=\{[B_{1},B_{2}]=0\}$.
Recall also that $\tX_{2}$ consists of four orbits $G.x^{(A)}$, $A=\I,\II,\III,\IV$ (\cite[Corollary 8.9]{Choy1}).
We notice that $\tX_{0}\setminus\tX_{1}=G.x\, (=G.x^{(\I)})$ because $\tX_0$ does not contain the other orbits $G.x^{(A)}$, $A=\II,\III,\IV$ in $\tX_{2}$ as
    $$
    B_1^{(\II)}=B_1^{(\III)}=B_1^{(\IV)}
    =\begin{pmatrix} 0 & 0 & 0 & 1 \\
    0&  0 & 0 & 0
    \\
    0 & -1 & 0 & 0
    \\
    0 & 0 & 0 &0
    \end{pmatrix}
    $$
are nilpotent matrices (hence with the only eigenvalue $0$).
By the dimension reason, $G.x$ is Zariski open dense in $\tX_{0}$.


\subsection{Factorization property}

In this subsection we introduce factorization property \cite[Section 2.16]{BFG} as a preliminary step of the proof of Proposition \ref{prop: smooth morphism}.
This property can be stated more generally, so we assume $V=\cc^k,W=\cc^N$ (symplectic and orthogonal vector spaces respectively).
The argument here, taken from \cite[Section 2.6]{Choy2}, originates from the proof of \cite[Lemma 12.3.2]{MO}.

For any given partition $\eta=(\eta_{1},\eta_{2},...,\eta_{e})$ of $k$, we denote by $S^{\eta}\cc$ the product of symmetric products $S^{\eta_{n}}\cc$.
Let $(S^{\eta}\cc)_{0}$ be the Zariski open subset of $S^{\eta}\cc$ consisting of the $e$-tuples $(Z_{1},Z_{2},...,Z_{e})$ such that the supports $\Supp(Z_{n})$ are mutually disjoint.
Let $\Delta S^{\eta}\cc$ be the Zariski closed subset of $S^{\eta}\cc$ consisting of $(Z_{1},Z_{2},...,Z_{e})$ such that $\Supp(Z_{n})$ are one-point sets for all $n$.
Let $(\Delta S^{\eta}\cc)_{0}:=(S^{\eta}\cc)_{0}\cap \Delta S^{\eta}\cc$.
Let $A_{\eta}\colon S^{\eta}\cc\to S^{k}\cc$, $(Z_{1},Z_{2},...,Z_{e})\mapsto \sum_{n=1}^{e}Z_{n}$.
It is a quasi-finite map and the restriction $A_{\eta}|_{(S^{\eta}\cc)_{0}}$ is \'etale.
There is a decomposition of $S^{k}\cc$ of the form
	\begin{equation}\label{eq: decomposition of Skcc}
	S^{k}\cc=\bigsqcup_{|\eta|=k} A_{\eta}((\Delta S^{\eta}\cc)_{0}).
	\end{equation}
Here $|\eta|=\sum_{n}\eta_{n}$.

Let $\pi\colon\bN_{V,W}\to S^{k}\cc$, $(B_{1},B_{2},i)\mapsto E_k(B_{1})$.
Using $\pi,A_{\eta}$, we consider fibre products $\bN_{V,W}\times_{S^{k}\cc}(S^{\eta}\cc)_{0}$ and $ \mu\inv(0)\times_{S^{k}\cc} (S^{\eta}\cc)_{0}$.
We notice that unless all $\eta_{n}$ are even, these fibre products are empty, since any generalized eigenspace of $B_1$ is a symplectic subspace (\cite[Lemma 5.2 (1)]{Choy1}).

From now on we may assume $\eta_{n}$ are all even.
We fix mutually orthogonal symplectic subspaces $V_{1}^{0},V_{2}^{0},...,V_{e}^{0}$ in $V$ with $\dim V_n^0=\eta_n$.
Let $\mu_{n}$ be the moment map with respect to $\Sp(V_{n}^{0})$ on $\bN_{V_{n}^{0},W}$ for each $n$.
Let
	\begin{equation}\label{eq: X'Y'}	
 	X':=\left(\bN_{V_{1}^{0},W}\oplus\bN_{V_{2}^{0},W}\oplus\cdots
	\oplus\bN_{V_{e}^{0},W}\right)\times_{S^{\eta}\cc}(S^{\eta}\cc)_{0},
	\quad
	Y':=\bN_{V,W}\times_{S^{k}\cc}(S^{\eta}\cc)_{0}
	\end{equation}
for short.
Note that $X'$ is Zariski open in $\bN_{V_{1}^{0},W}\oplus\bN_{V_{2}^{0},W}\oplus\cdots\oplus\bN_{V_{e}^{0},W}$.
Using the orthogonal decomposition $V=V_{1}^{0}\oplus V_{2}^{0}\oplus\cdots\oplus V_{e}^{0}$, we define a morphism
	\begin{equation}\label{eq: pi'eta}
	\textstyle
	\pi_{\eta}'\colon X'\to Y',\quad
 	(B_{n,1},B_{n,2},i_{n})_{1\le n\le e}
	\mapsto (\bigoplus_{n=1}^{e}B_{n,1},B_{2},\bigoplus_{n=1}^{e}i_{n},Z_{1},Z_{2},...,Z_{e}),
	\end{equation}
where the matrix $B_{2}$ is determined so that its block $B_{2}^{(m,n)}\colon V_{n}^{0}\to V_{m}^{0}$ satisfies
 	\begin{equation}\label{eq: factorization}
	\left\{
	\begin{array}{lll}
	B_{2}^{(n,n)}=  B_{n,2} &&
	\\
	B_{m,1}B_{2}^{(m,n)}-B_{2}^{(m,n)}B_{n,1}+i_{m}i_{n}^{*}=0
	&& \mbox{if $m\neq n$}
	\end{array}	
	\right.
	\end{equation}
and $Z_{n}$ is the set of eigenvalues of $B_{1}|_{V_{n}^{0}}$ counted with multiplicities.
Note that $B_{2}^{(m,n)}$ $(m\neq n)$ is determined uniquely by \eqref{eq: factorization} as $B_{m,1}$ and $B_{n,1}$ have different eigenvalues.
Hence each matrix element of $B_2$ is given as a rational function in the matrix elements of $B_{n,1},B_{n,2},i_{n}$ which is defined on $X'$.
Note that $B_2\in \fp(V)$.
This can be seen by taking the right adjoint to each equation in  \eqref{eq: factorization}.
Thus we deduce the following lemma:

\begin{lem}\label{lem: cl embed}
$\pi_\eta'$ is a closed embedding.
Moreover $(\pi_\eta')^*\colon \cc[Y']\to \cc[X']$ is a surjective algebra homomorphism with the kernel ideal generated by the second equation of \eqref{eq: factorization}, where $\cc[X'],\cc[Y']$ denote the affine coordinate rings of $X',Y'$ respectively.
\end{lem}

\proof
Note that $X',Y'$ are affine varieties since $(S^\eta\cc)_0$ and $A_\eta((S^\eta\cc)_0)$ are affine varieties and $A_\eta|_{(S^\eta\cc)_0}$ is a finite morphism.
The lemma itself is a direct consequence of the definition of $\pi_\eta'$.
\qed\vskip.3cm

We define closed subschemes of $X',Y'$ as the $0$ loci of the moment maps:
	$$
	 X:=\left(\mu_{1}\inv(0)\times\mu_{2}\inv(0)\times\cdots\times\mu_{e}\inv(0)\right)\times_{S^{\eta}\cc}(S^{\eta}\cc)_{0},
	\quad
	Y:=\mu\inv(0)\times_{S^{k}\cc}(S^{\eta}\cc)_{0}.
	$$
Since $\mu$ becomes the product $\prod_{n=1}^{e}\mu_{n}$ on $X'$ after composition with $\pi_{\eta}'$ and the \'etale projection $Y'\to\bN_{V,W}$, we have $\pi'_{\eta}(X)\subset Y$.
We denote this restriction by $\pi_{\eta}\colon X\to Y$.

Let $\Sp:=\Sp(V),\Sp':=\prod_{n=1}^{e}\Sp(V_{n}^{0})$ with the Lie algebras $\fsp,\fsp'$ respectively.
The $\Sp$-action on $\bN_{V,W}$ naturally lifts to $Y'$ (by giving the trivial $G$-action on the factor $S^{\eta}\cc$).
Similarly there is also a natural $\Sp'$-action on $X'$.
We consider another $\Sp'$-action on $\Sp\times X'$, defined by $h.(g,x)=(gh\inv,h.x)$.

\begin{lem}\label{lem: quotient}
$(1)$
$\pi_\eta,\pi_\eta'$ are $\Sp'$-equivariant closed embeddings.

$(2)$
The action map $\sigma\colon \Sp\times X\to Y,\ (g,x)\mapsto g.\pi_\eta(x)$ is a surjective smooth morphism with fibres isomorphic to $\Sp'$.
\end{lem}

\proof
(1)
By Lemma \ref{lem: cl embed} we know that $\pi_\eta'$ is a closed embedding with ideal defined by \eqref{eq: factorization}.
The restriction to $\mu=0$ gives further equations on $B_{1}$ and the diagonal blocks $B_{n,1}$.
Taking quotients by the ideals corresponding to $(\pi_\eta')^*$ we obtain the induced surjective homomorphism $(\pi_\eta)^*\colon \cc[Y]\to \cc[X]$, so that $\pi_\eta$ is a closed embedding.

We check that $\pi_\eta,\pi_\eta'$ are $\Sp'$-equivariant.
Let $h=(h_1,h_2,...,h_e)\in \Sp'$.
We need to check the block corresponding to $\bigoplus_{n=1}^eh_n.B_{n,2}$ is $h_mB_2^{(m,n)}h_n\inv$.
This follows from the uniqueness of the solution of \eqref{eq: factorization} after conjugating the first equation \eqref{eq: factorization} by $h_n$ and multiplying the second equation \eqref{eq: factorization} by $h_m$ on the left and $h_n\inv$ on the right.
The details are left to the reader (or see the proof of \cite[Lemma 2.7]{Choy2}).

(2)
We denote the $\Sp$-action map by $\sigma'\colon \Sp\times X' \to Y'$, $(g,x)\mapsto g.\pi_{\eta}'(x)$.
We claim that $\sigma'$ is a submersion with fibres isomorphic to $\Sp'$.
Let us finish the proof of the lemma assuming this claim.
We notice that $Y$ is a closed subscheme of the image variety of $\sigma'$.
For, $\Image(\sigma')$ is a subscheme defined only by the second equation of \eqref{eq: factorization} while $Y$ is defined by the additional equations from $\mu=0$ on the diagonal blocks.
Thus $\sigma$ satisfies the statement as it comes from $\sigma'$ via a base change.

We now prove the claim in the rest of the proof.
First we describe the fibres of $\sigma'$.
For the second $\Sp'$-action on $\Sp\times X'$, let $\bO_{(g,x)}$ be the free $\Sp'$-orbit of $(g,x)$ in $\Sp\times X'$.
Note that set-theoretically
	\begin{equation}\label{eq: fibre}
	\sigma'{}\inv(\sigma'(g,x))=\bO_{(g,x)}.
	\end{equation}
For, this amounts to $\sigma'(g,x)=g.\pi_{\eta}'(x)\in \pi_{\eta}'(X')$ implies $g\in \Sp'$.
Let us prove this.
We write $x=(B_{1},B_{2},i,Z_{1},Z_{2},...,Z_{e})\in \pi_\eta'(X')$.
The sum of generalized eigenspaces of $B_{1}$ with eigenvalues in $Z_{n}$ is $V_{n}^{0}$ for each $n=1,2,...,e$.
The assumption $g.\pi_{\eta}'(x)\in \pi_{\eta}'(X')$ asserts that $g.B_{1}$ also has $V_{n}^{0}$ as the sum of generalized eigenspaces with eigenvalues in $Z_{n}$ for each $n$.
This implies $g\inv(V_{n}^{0})=V_{n}^{0}$ and thus we have $g\in \Sp'$.

It remains to prove that $\sigma'$ is submersive: the differential $d\sigma'$ at any point has constant rank $\dim \Sp+\dim X'-\dim \Sp'$.
Equivalently $\Ker(d\sigma'{}_{(g,x)})=T_{(g,x)}\bO_{(g,x)}$.
By the $\Sp$-equivariancy we may assume $g$ is the identity element of $\Sp$.
A direct computation shows
    $$
    \Ker (d\sigma'_{(g,x)})=\left\{(\xi,\delta x)\in \fsp\oplus \bigoplus_{n=1}^e\bN_{V^0_n,W} \middle|\xi.x+\delta x=0\right\}=T_{(g,x)}\bO_{(g,x)}.
    $$
This completes the proof.
\qed\vskip.3cm

If $\eta=(k)$, we have $X'=Y', X=Y$ and $\Sp=\Sp'$, so the above lemma gives only the obvious isomorphism.
Otherwise the \'etale open subset of $\mu\inv(0)$ associated to $\eta$ is isomorphic to a product of strata involving strictly lower instanton numbers.
We call this property the \textit{factorization property} of $\mu\inv(0)$.
\vskip.3cm


\subsection{Proof of Proposition \ref{prop: smooth morphism}}\label{subsec: proof of prop}
Let $k=4,N=3$ and $\eta=(2,2)$.
Let $q\colon Y\to \mu\inv(0)$ be the (open \'etale) projection.
We set the morphism in Proposition \ref{prop: smooth morphism} (1) to be the restriction of $q\sigma$ to $(q\sigma)\inv(\tX_{0})$.
Due to Lemma \ref{lem: quotient} (2), we need to show that $(q\sigma)\inv(\tX_{0})\cong\Sp\times\cc^{*}\times\cc\times(\rho\inv(0))^{2}$.
By \cite[Theorem 4.1 (1)]{Choy1}, there is a natural identification $\mu_{2}\inv(0)=\cc'\times\cc\times\rho\inv(0)$, where $\cc'=\{\tr(B_{1})\}\, (=\cc)$.
By the factorization property, $(q\sigma)\inv(\tX_{0})$ is the closed subscheme of $X=(\mu_{2}\inv(0))^{2}\times_{S^{\eta}\cc}(S^{\eta}\cc)_{0}$ whose $(\cc')^{2}$- and $\cc^{2}$-factors both consist of points of the form $(a,-a)$, $a\in\cc$.
This proves the first item.

The second item follows from $G.x=\tX_{0}\setminus\tX_{1}$, so that $\rho(i)\neq0$ whenever $i$ comes from $G.x$.
This finishes the proof of Proposition \ref{prop: smooth morphism}.
 

\subsection{Remark}
\label{subsec: remark}
Note that $B^{(\I)}=v_5$ has the $1/2$-eigenspace $\cc\langle e_1+e_3,e_2+e_4\rangle$ and the $-1/2$-eigenspace $\cc\langle e_1-e_3,e_2-e_4\rangle$.
Thus by choosing any $g\in \Sp$ which maps these $\pm1/2$-eigenspaces to $\cc\langle e_1,e_2\rangle,\cc\langle e_3,e_4\rangle$ respectively, the ADHM datum
$g.(B^{(\I)}_1,B^{(\I)}_2, i^{(\I)})$ comes from the product $(\mu_{(2)}\inv(0))^2$ by the factorization property.

The stabilizer in $\Sp(V)$ of any $x=(B_{1,1},B_{1,2},i_1)\oplus (B_{2,1},B_{2,2},i_2) \in(\mu_{(2)}\inv(0))^2\times_{S^{(2,2)}\cc}(S^{(2,2)}\cc)_0$ is easy to compute, because $\fp(\cc^2)=\cc$, and thus $\Sp(V)^{B_{1,1}\oplus B_{2,1}}=\Sp(1)^2$.
Therefore $\Sp(V)^x=\Sp(1)^{i_1}\times\Sp(1)^{i_2}$.
Since $\rank(i_{1})=\rank(i_{2})=1$ (as $i_{1},i_{2}\in\rho\inv(0)$), both $\Sp^{i_{1}},\Sp^{i_{2}}$ are isomorphic to the additive group $\cc$.
This is what we computed in Section \ref{subsec: corrected lemma} at the Lie algebra level.


\section{Scheme structure of the GIT quotient $\tX_{0}\git\Sp$}
\label{sec: scheme structure of GIT quotient}

Recall that $\tX_0$ is the complement of the vanishing locus of a $G$-invariant polynomial, the discriminant of the characteristic polynomial of $B_1$.
So it is an affine open subset of $\tX$ and the GIT quotient $\tX_0\git \Sp$ is also an affine open subset of $\tX\git \Sp$.

In this section we further describe the scheme structure of $\tX_0\git \Sp$.
\begin{lem}\label{lem: GIT}
The GIT quotient $\rho\inv(0)\git \Sp$ is a non-reduced scheme.
\end{lem}

\begin{cor}\label{cor: GIT}
$\tX_{0}\git \Sp$ is a non-reduced scheme.
\end{cor}


\subsection{Proof of Corollary \ref{cor: GIT}}
We first prove the corollary assuming Lemma \ref{lem: GIT}.

By Proposition \ref{prop: smooth morphism} we have an isomorphism $\tX_{0}\cong (\Sp\times \cc^*\times \cc\times (\rho\inv(0))^2)/(\Sp'\times\zz/2\zz)$.
Here the $\Sp'$-action was introduced before Lemma \ref{lem: quotient}, so $\Sp'$ acts nontrivially only on the two factors $\Sp,(\rho\inv(0))^{2}$.
And $\zz/2\zz$ acts only on the factor $\cc^{*}$ by $a\mapsto -a$.
These two actions commute and the $\Sp'\times \zz/2\zz$-action is free.
Considering the quotient by the further $\Sp$-action action on the factor $\Sp$, we obtain an isomorphism $\tX_0\git \Sp\cong (\Sp\times \cc^*\times \cc\times (\rho\inv(0))^2)\git (\Sp\times \Sp'\times\zz/2\zz)$.
Taking the quotients by $\Sp$ and then by $\Sp'\times\zz/2\zz$, we obtain an isomorphism
	$$
	\tX_0\git \Sp\cong \cc^*\times \cc\times(\rho\inv(0)\git \Sp(1))^2.
	$$
Since the factor $\rho\inv(0)\git \Sp(1)$ is non-reduced by Lemma \ref{lem: GIT}, our assertion is proven.


\subsection{Proof of Lemma \ref{lem: GIT}}

Recall that from the $ab$-diagram of elements in $\rho\inv(0)$ the reduced scheme $\rho\inv(0)\git\Sp(1)_{\mathrm{red}}$ is a one-point set (\cite[Theorem 4.1 (5)]{Choy2}).

Let us check that $\rho\inv(0)\git\Sp(1)$ is a non-reduced subscheme of $W\dual$.
On $\Hom(W,V)=W\dual\otimes V$, the defining equations $\rho(i)=ii^*=0$ are quadratic in $W\dual$.
Therefore the $\Sp$-invariants in the ideal always have degree $\ge2$, as the invariants are obtained from the integration over $\Sp$ with respect to the normalized Haar measure (i.e.\ Reynolds' operator).
So the invariant ideal in $\cc[W\dual]$ is contained in $\fm_{0}^2$, where $\fm_{0}$ denotes the maximal ideal defining $0\in W\dual$.
This means that $\rho\inv(0)\git \Sp(1)$ contains the fat point $\cc[x,y,z]/(x,y,z)^2$ as a closed subscheme, hence it is non-reduced.


\subsection{Remark}

In fact $\rho\inv(0)\git \Sp(1)$ coincides with the above fat point.
This can be checked by applying the integration over $\Sp$ to the matrix elements of $ii^{*}$.
We do not give the detailed computation here (see \cite[Appendix A]{Choy3}).